\documentclass[11pt,twoside]{article}
\usepackage{latexsym,amssymb,amsmath}
\usepackage{amsfonts}
\RequirePackage{graphics,epsfig,psfrag}%,picins}
\def\abs#1{\left \vert #1 \right \vert}
\def\hat#1{\widehat{#1}}
\def\RR{{\bf R}} %reelle Zahlen
\def\ZZ{{\bf Z}} %reelle Zahlen
\def\QQ{{\bf Q}} %reelle Zahlen
\def\SS{{\bf S}} %reelle Zahlen
\def\GG{{\mathbf{G}}} %groupe G
\def\Ga{{\mathbf{\Gamma}}} %groupe Gamma
\def\cS{{\cal S}} %reelle Zahlen
\def\Mod#1{\,(\hbox{\rm mod}\,#1)}
\def\cn{\hbox{\rm cn}\,}
\def\H{{\rm H}}
\def\id{{\mathbf{Id}}}

\def\phi{\varphi}
\def\eps{\varepsilon}

\def\cC{{\cal C}}

\def\cT{{\cal T}}
\def\pn{\medskip\par\noindent}
\def\bi{\vspace{-2pt}\begin{itemize}\itemsep -2pt plus 1pt minus 1pt}
\def\ei{\end{itemize}\vspace{-4pt}}
\def\bn{\vspace{-2pt}\begin{enumerate}\itemsep -2pt plus 1pt minus 1pt}
\def\en{\end{enumerate}\vspace{-4pt}}
\newcommand{\Pf}{{\em Proof}. }

\newcommand{\EPf}{\hbox{}\hfill$\Box$\vspace{.5cm}}
\def\[#1\]{\begin{eqnarray}#1\end{eqnarray}}
\def\$#1\${\begin{eqnarray*}#1\end{eqnarray*}}
\def\p{\vspace{-5.17ex} \hspace{6.5ex}}

\def\pent#1#2{\lfloor\frac{#1}{#2}\rfloor}

\def\comb#1#2{{}{\textstyle{{#1} \choose {#2}}}}
\def\sign#1{{\rm sign}\,\bigl( #1 \bigr)}
\def\Sum{\mathop{\sum}\limits}

\def\abs#1{\left \vert #1 \right \vert}

\def\frac#1#2{{\textstyle{{#1} \overwithdelims.. {#2}}}}
\def\Frac#1#2{{\displaystyle{{#1} \overwithdelims.. {#2}}}}

\catcode`\@=11
\def\system#1{\left\{\null\,\vcenter{\openup\jot\m@th
\ialign{
\strut\hfil$\displaystyle{##}$&
$\,\displaystyle{{}##}\,$\hfil&&
\strut\hfil$\,\displaystyle{##}$&
$\,\displaystyle{{}##}\,$\hfill
\hfil\crcr#1\crcr}}\right.}

\def\cmatrix#1{\left [
\null\,\vcenter{%\openup\jot\m@th
\ialign{
\hfil${##}\ $\hfil &
\hfil$\ {##}\ $\hfil&&
\hfil$\ {##}\ $\hfil&
\hfil$\ {##}$\hfil
\crcr#1\crcr}}\right ]}

\def\@opargbegintheorem#1#2#3{\par\addvspace{6pt plus3pt minus2pt}%
    \def\@tempa{#3}%
    \noindent{\bf #1 #2 \ifx\@tempa\empty\unskip\else\unskip\ (#3).\fi\hskip.5em}\csname#1font\endcsname\ignorespaces
\ignorespaces}

\def\@endtheorem{\par\addvspace{6pt plus3pt minus2pt}}
\def\@begintheorem#1#2#3{\par\addvspace{8pt plus3pt minus2pt}%
              \noindent{\csname#1headfont\endcsname#1\ \ignorespaces#3 #2.}%
              \csname#1font\endcsname\hskip6pt\ignorespaces}
\def\@endtheorem{\par\addvspace{8pt plus3pt minus2pt}\@endparenv}

\catcode`\@=12

\newtheorem{theorem}{Theorem}[section]
\newtheorem{thm*}{Theorem}
\newtheorem{corollary}[theorem]{Corollary}
\newtheorem{lemma}[theorem]{Lemma}
\newtheorem{proposition}[theorem]{Proposition}
\newtheorem{definition}[theorem]{Definition}
\newtheorem{remark}[theorem]{Remark}
\newtheorem{example}[theorem]{Example}
\newtheorem{examples}[theorem]{Examples}

\date{\today}
\begin{document}
\pagestyle{myheadings}
\markboth{P. -V. Koseleff, D. Pecker}{{\em Chebyshev diagrams for rational knots}}
%%%%%%%%%%%%%%%%%%%%% Publisher's Area please ignore %%%%%%%%%%%%%%
%wis%\catchline{}{}{}{}{}
%%%%%%%%%%%%%%%%%%%%%%%%%%%%%%%%%%%%%%%%%%%%%%%%%%%%%%%%%%%%%%%%%%%
\title{Chebyshev diagrams for rational  knots}
\author{P. -V. Koseleff, D. Pecker}
\maketitle
\begin{abstract}
We show that every rational knot $K$ of crossing number $N$ admits
a polynomial parametrization
$x=T_a(t), \, y = T_b(t), z = C(t)$ where $T_k(t)$ are the
Chebyshev polynomials, $a=3$ and $b+ \deg C = 3N.$ We show that
every rational knot also admits  a polynomial parametrization with $a=4$.
If $C (t)= T_c(t)$ is a Chebyshev polynomial, we call such a knot a harmonic
knot. We give the classification of  harmonic knots for $a \le 4.$
%Most results are derived from continued fractions and their matrix representations.
\pn {\bf keywords:}{
Polynomial curves, %Chebyshev polynomials,
Chebyshev curves, rational knots, continued fractions}
\end{abstract}
\begin{center}
\parbox{12cm}{\small
\tableofcontents
}
\end{center}
\vspace{1cm}
\section{Introduction}
We study the polynomial parametrization of knots,
viewed as non singular space curves.
Vassiliev proved
that any  knot  can be represented by a polynomial
embedding $\RR \to \RR ^3 \subset {\bf S}_3 $ (\cite{Va}).
Shastri (\cite{Sh}) gave another proof of this theorem,
he also found explicit  parametrizations of the trefoil
and of the figure-eight knot (see also \cite{Mi}).
\pn
We shall study polynomial embeddings of the  form
$x=T_a(t), \, y= T_b(t), \, z=C(t)$ where $a$ and $b$ are coprime
integers and $T_n(t)$ are  the classical Chebyshev polynomials defined by
$T_n(\cos t) = \cos nt$.
The projection of such a curve on the $xy$-plane is
the Chebyshev  curve  $\cC(a,b): T_b(x)=T_a(y)$ which has
exactly $\frac 12 (a-1)(b-1)$ crossing points (\cite{Fi,P1,P2}).
We will say that such a knot has the Chebyshev diagram $\cC(a,b)$.
\pn
We observed in \cite{KP1} that  the trefoil can be parametrized by
Chebyshev polynomials: $x=T_3(t);\, y=T_4(t);\,z= T_5(t)$.
This led us to study Chebyshev knots in \cite{KP3}.
\begin{definition}
A knot in $\RR ^3  \subset  {\SS}^3$ is the Chebyshev knot $\cC(a,b,c,\phi)$
if it admits the one-to-one parametrization
$$ x=T_a(t); \  y=T_b(t) ; \    z=T_c(t + \phi) $$
where $t \in \RR$, $a$ and $b$ are coprime integers, $c$
is an integer and $\phi$ is a real constant.
\smallskip\par\noindent
When $\phi=0$ and $a,b,c$ are coprime, it is
denoted by $\H(a,b,c)$ and is called a harmonic knot.
\end{definition}
We proved that any knot is a Chebyshev knot.
Our proof uses theorems on braids by Hoste, Zirbel and Lamm (\cite{HZ}),
and a density argument.
In a joint work with F. Rouillier (\cite{KPR}), we developed an effective method
to enumerate all the knots $\cC(a,b,c,\phi), \phi \in \RR$ where $a=3$ or $a=4$, $a$ and
$b$ coprime.
\pn
Chebyshev knots are polynomial analogues  of
Lissajous knots that admit a  parametrization of the form
$$
x=\cos (at ); \   y=\cos (bt + \phi) ; \  z=\cos (ct + \psi)
$$
where $ 0 \le t \le 2 \pi $ and where $ a, b, c$ are pairwise
coprime integers. These knots, introduced  in \cite{BHJS},
have been studied by  V. F. R. Jones, J. Przytycki,
C. Lamm, J. Hoste and L. Zirbel. Most known properties of Lissajous
 knots are deduced from their symmetries
(see \cite{BDHZ,Cr,HZ,JP,La1}).
\pn
The symmetries of harmonic knots, obvious from the parity of Chebyshev polynomials,
are different from those of Lissajous.
For example, the figure-eight knot which is amphicheiral but not a Lissajous knot, is
the harmonic knot $ \H(3,5,7).$
\pn
We proved that the harmonic knot $\H(a,b, ab-a-b)$ is alternate, and deduced
that there  are infinitely many amphicheiral harmonic knots and
infinitely many strongly invertible harmonic knots.
We also proved (\cite{KP3}) that the torus knot $T(2, 2n+1)$ is the harmonic knot
$ \H(3,3n+2,3n+1)$.
\pn
In this article, we give the classification of the harmonic knots $\H(a,b,c)$ for $a \le 4.$
We also  give explicit polynomial parametrizations  of all rational knots.
The diagrams of our knots are Chebyshev  curves of minimal degrees with a small number
of crossing points. The degrees of the
height polynomials are small.
\pn
In section {\bf \ref{cf}.} we recall the Conway notation for  rational knots,
and the computation of their Schubert fractions  with continued fractions.
We observe that Chebyshev diagrams  correspond
to continued fractions of the form  $[\pm 1, \ldots, \pm 1]$ when $a=3$
and of the form  $[\pm 1,\pm 2,  \ldots, \pm 1, \pm 2]$ when $a=4$.
We show results on our continued fraction expansion:
\pn
{\bf Theorem \ref{th1}.}\\
{\em Every rational number $r$ has a unique continued fraction expansion
$r = [e_1, e_2, \ldots, e_n ]$, $e_i= \pm 1$,
where there are no two consecutive sign changes in the sequence $(e_1,\ldots, e_n ).$ }
\pn
We have a similar theorem for continued fractions of the form
$r = [\pm 1, \pm 2, \ldots, \pm 1, \pm 2]$.
We provide a formula for the crossing number of the corresponding knots.
Then we study the matrix interpretation of these continued fraction expansions.
As an application, we give optimal Chebyshev diagrams for the torus knots
$T(2,N)$, the twist knots $\cT_n$, the generalized stevedore knots and some
others.
\pn
In section {\bf \ref{harmonic}.} we describe the harmonic knots $\H(a,b,c)$ where
$a \le 4.$
We begin with a careful analysis of the nature of the crossing points,
giving  the Schubert fractions of $\H ( 3,b,c) $ and $\H(4,b,c).$
Being rather long,  the proofs of these results will be given in the last paragraph.
We deduce the following algorithmic classification theorems.
\pn
{\bf Theorem \ref{h3bc}.}\\
{\em Let $K= \H(3,b,c)$.
There exists a unique pair $(b',c')$ such that (up to mirror symmetry)
$$ K = \H (3,b',c'), \   b'<c'< 2b', \   b' \not\equiv c' \Mod 3.$$
The crossing number of $K$ is $ \Frac 13 ( b'+c')$.\\
The Schubert fractions $ \Frac \alpha \beta $ of $K$ are such that
$\beta ^2 \equiv \pm 1 \Mod \alpha.$
}
\pn
{\bf Theorem \ref{h4bc}.}\\
{\em
Let  $K= \H(4,b,c).$
There exists a unique pair $(b',c') $ such that (up to mirror symmetry)
$$ K= \H(4,b',c'), \  b' < c' < 3b', \ b' \not\equiv c' \Mod 4.$$
The crossing number of $K$ is $ \Frac 14 ( 3b'+c'-2)$.\\
$K$ has a Schubert fraction $ \Frac \alpha \beta $ such that
$ \beta ^2 \equiv \pm 2 \Mod \alpha.$
}
\pn
We notice that the trefoil is the only knot which is both of form
$\H(3,b,c)$ and $\H(4,b,c).$
We remark that the $6_1$ knot (the stevedore knot) is not a harmonic knot
$\H( a,b,c), \  a \le 4.$
\pn
In section {\bf \ref{diagrams}.} we find explicit polynomial parametrizations
of all rational knots. We first compute the optimal Chebyshev diagrams
for $a=3$ and $a=4$. Then we define a height polynomial of small degree.
More precisely:
\pn
{\bf Theorem \ref{gauss3}.}\\
{\em Every rational knot of crossing number $N$ can be parametrized by
$x=T_3(t), y= T_b(t), z= C(t)$ where $b + \deg C = 3N$.
Furthermore,  when the knot
is amphicheiral, $b$ is odd and we can choose $C$ to be an odd polynomial.}
\pn
In the same way we show:
{\em Every rational knot of crossing number $N$ can be parametrized by
$x=T_4(t), y= T_b(t), z= C(t)$} where $b$ is odd and $C$ is an odd polynomial.
\pn
As a consequence, we see that any rational knot has a representation $K \subset \RR^3$ such that
$K$ is symmetrical about the $y$-axis (with reversed orientation).
It clearly implies the classical result:
{\em every rational knot is strongly invertible}.
\pn
We give  polynomial parametrizations of
the torus knots $T(2,2n+1).$
We also give the first polynomial parametrizations of
the  twist knots and the generalized stevedore knots.
We conjecture that the lexicographic degrees of our polynomials are minimal
(among odd or even polynomials).
%%%%%%%%%%%%%%%%%%%%%%%%%%%%%%%%%%%%%%%%%%%%%%%%%%%%%%%%%%%%%%%%%%%%%%%%%%%%%%%%%%%%%%%%%%%%%%%%%%%
\section{Continued fractions and rational Chebyshev knots}\label{cf}
\begin{figure}[th]
\begin{center}
\begin{tabular}{cccc}
{\scalebox{.6}{\includegraphics{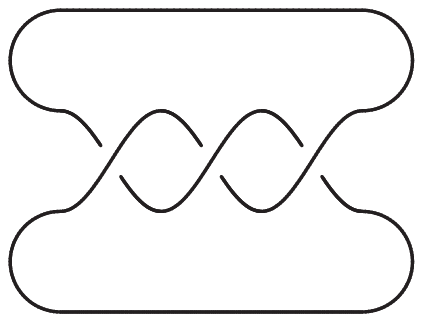}}}&
{\scalebox{.6}{\includegraphics{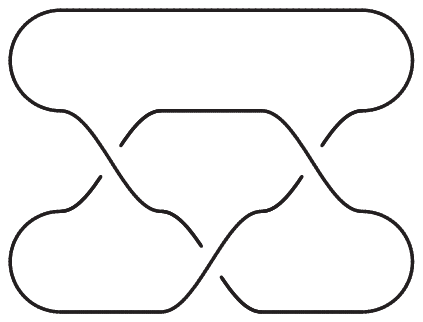}}}&
{\scalebox{.6}{\includegraphics{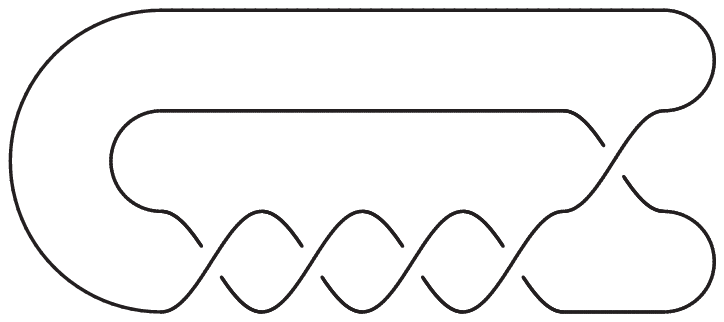}}}&
{\scalebox{.6}{\includegraphics{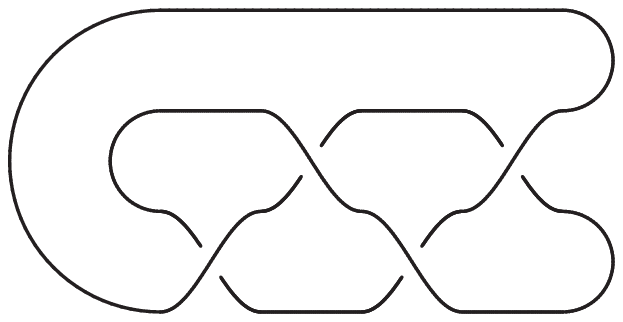}}}\\
$C(3)$&$C(-1,-1,-1)$&$C(4,-1)$&$C(1,1,-1,-1)$
\end{tabular}
\end{center}
\caption{Diagrams of the standard trefoil}\label{trefoils}
\end{figure}

%\subsection{Two-bridge knots}
A two-bridge knot (or link) admits a diagram in Conway's normal form.
This form, denoted by
$C(a_1, a_2, \ldots, a_n)$  where $a_i$ are integers, is explained by the following
picture (see \cite{Con}, \cite{Mu} p. 187).
\psfrag{a}{\small $a_1$}\psfrag{b}{\small $a_2$}%
\psfrag{c}{\small $a_{n-1}$}\psfrag{d}{\small $a_{n}$}%
\begin{figure}[th]
\begin{center}
{\scalebox{.8}{\includegraphics{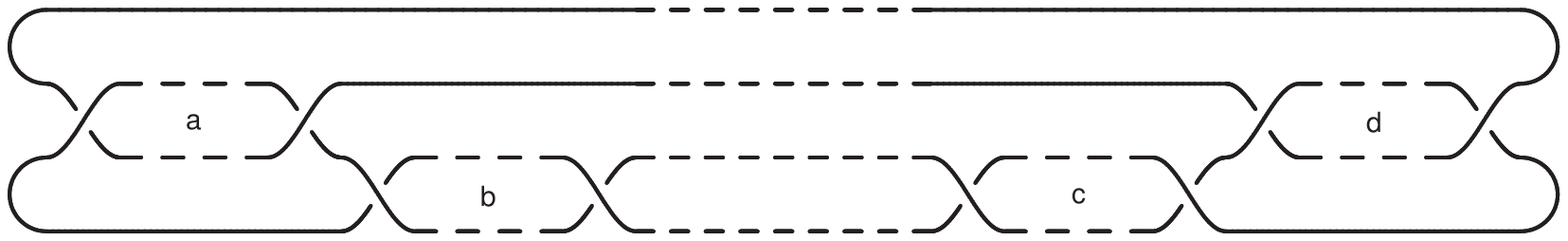}}}\\[30pt]
{\scalebox{.8}{\includegraphics{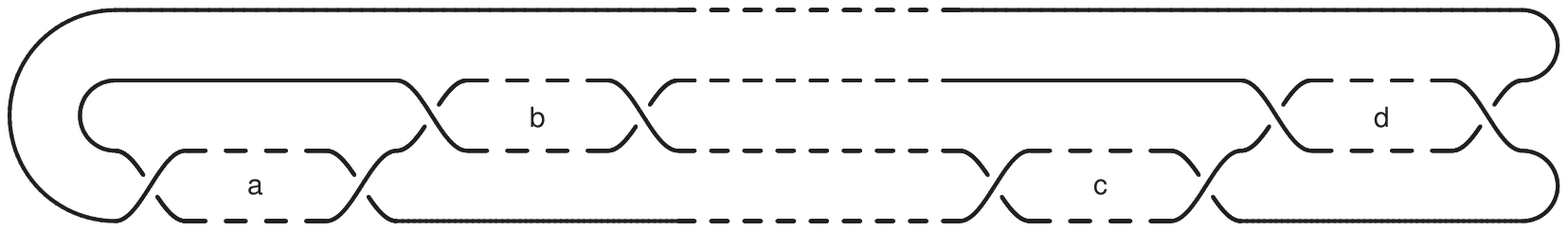}}}
\end{center}
\caption{Conway's normal forms}
\label{conways3}
\end{figure}
The number of twists is denoted by the integer
$\abs{a_i}$, and the sign of $a_i$ is defined
as follows: if $i$ is odd, then the right twist is positive,
if $i$ is even, then the right twist is negative.
On  Fig. \ref{conways3} the $a_i$ are positive (the $a_1$ first twists are right twists).
\pn
\begin{examples}
The trefoil has the following  Conway's normal forms
$C(3)$, $C(-1,-1,-1)$, $C(4, -1)$ and $C(1,1,-1,-1 ).$
The diagrams in Figure \ref{trefoils} clearly represent the same trefoil.
\end{examples}
\def\p{{\small $+$}}
\def\m{{\small $-$}}
\pn
The two-bridge links are classified by their Schubert fractions
$$
 \Frac {\alpha}{\beta} =
a_1 + \Frac{1} {a_2 + \Frac {1} {a_3 + \Frac{1} {\cdots +\Frac 1{a_n}}}}=
[ a_1, \ldots, a_n], \quad \alpha >0.
$$
We shall denote  $S \bigl( \Frac {\alpha}{\beta} \bigr)$  a two-bridge link with
Schubert fraction $ \Frac {\alpha}{\beta}.$
The two-bridge  links
$ S (\Frac {\alpha} {\beta} )$ and $ S( \Frac {\alpha ' }{\beta '} )$ are equivalent
if and only if $ \alpha = \alpha' $ and $ \beta' \equiv \beta ^{\pm 1} ( {\rm mod}  \  \alpha).$
The integer $ \alpha$ is odd for a knot, and even for a two-component link.
If $K= S (\Frac {\alpha}{\beta} ),$ its mirror image is
$ \overline{K}= S ( \Frac {\alpha}{- \beta} ).$

We shall study knots with a Chebyshev diagram $\cC (3,b) : \  x= T_3(t), y= T_b(t).$
It is remarkable that such a diagram is already in Conway normal form
(see Figure \ref {conways3}).
Consequently, the Schubert fraction of such a knot is given by a
continued fraction of the form
$ [ \pm 1, \pm 1, \ldots,\pm 1 ].$
For example the only  diagrams of Figure \ref{trefoils}
which may be Chebyshev are the second and the last
(in fact they are Chebyshev).
\pn
Figure \ref{t7} shows  a typical example of  a knot with a Chebyshev diagram.
\def\p{{$\scriptstyle{+}$}}%
\def\m{{\small $-$}}%
\begin{figure}[th]
\begin{center}
\begin{tabular}{cc}
{\scalebox{.8}{\includegraphics{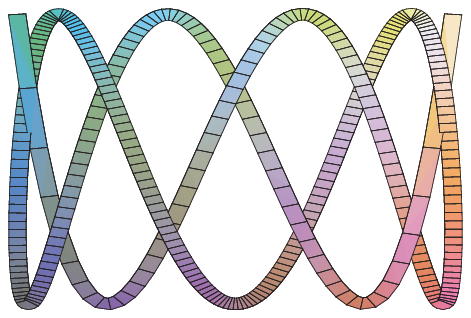}}}&
\psfrag{n1}{\m}\psfrag{n2}{\m}\psfrag{n3}{\m}%
\psfrag{n4}{\p}\psfrag{n5}{\p}\psfrag{n6}{\p}%
\psfrag{n7}{\m}\psfrag{n8}{\m}\psfrag{n9}{\m}%
{\scalebox{1}{\includegraphics{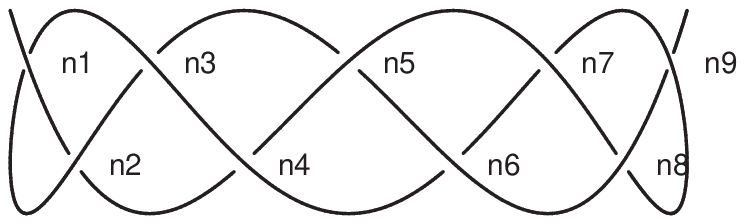}}}
\end{tabular}
\end{center}
\caption{A Chebyshev diagram of  the torus knot $T(2,7)$}
\label{t7}
\end{figure}
\pn
This knot is defined by $ x= T_3(t), \, y= T_{10}(t), \, z = -T_{11}(t).$
Its $xy$-projection is in the Conway normal form
$ C(-1,-1,-1, 1,1,1, -1,-1,-1).$ Its Schubert fraction is then $\Frac{7}{-6}$ and this knot is
the torus knot $T(2, 7)=S(\Frac{7}{-6}) = S(7)$.
\pn
We shall also need to study knots with a diagram illustrated by the
following picture.
\psfrag{a}{\small $b_1$}\psfrag{b}{\small $a_1$}\psfrag{c}{\small $c_1$}%
\psfrag{d}{\small $b_n$}\psfrag{e}{\small $a_n$}\psfrag{f}{\small $c_n$}%
\begin{figure}[th]
\begin{center}
{\scalebox{.8}{\includegraphics{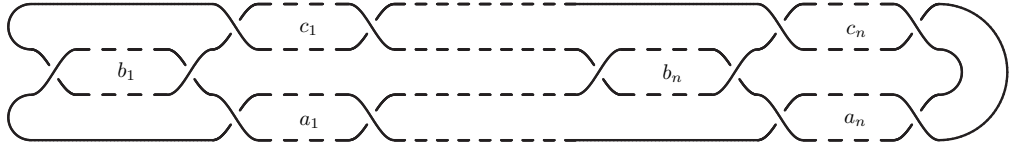}}}
\end{center}
\caption{A knot isotopic to $C(b_1,a_1+c_1,b_2,a_2+c_2,\ldots, b_n,a_n+c_n)$}
\label{conway4}
\end{figure}
In this case, the $a_i$ and the $c_i$ are positive if they are left twists,
the $b_i$ are positive if they are right twists (on our figure  $a_i, b_i, c_i $ are positive).
Such a knot is equivalent to a knot with Conway's normal form
$C(b_1, a_1+c_1, b_2, a_2 +c_2, \ldots, b_n, a_n+c_n )$
(see \cite {Mu} p. 183-184).
%, the integers $b_i$ may be  $0$ or $ \infty$.
\pn
We shall study  the knots with a Chebyshev diagram
$\cC(4,k): x= T_4(t), \, y= T_ k(t)$.
In this case we get diagrams of the form illustrated by Figure \ref{conway4}.
Consequently, such a knot  has a Schubert fraction of the form
$[ b_1, d_1, b_2, d_2,  \ldots, b_n, d_n ]$
with $b_i=\pm 1$,  $d_i = \pm 2$ or $d_i=0.$

\begin{figure}[th]
\begin{center}
\begin{tabular}{ccc}
\psfrag{n1}{\p}\psfrag{n2}{\m}\psfrag{n3}{\p}%
\psfrag{n4}{\p}\psfrag{n5}{\p}\psfrag{n6}{\p}%
{\scalebox{1.01}{\includegraphics{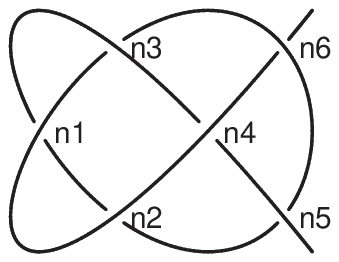}}}&
%\psfrag{n1}{\p}\psfrag{n2}{\p}\psfrag{n3}{\p}\psfrag{n4}{\p}%
%{\scalebox{.8}{\includegraphics{di4_1_3.eps}}}&
\psfrag{n1}{\m}\psfrag{n2}{\m}\psfrag{n3}{\m}%
\psfrag{n4}{\p}\psfrag{n5}{\p}\psfrag{n6}{\p}%
{\scalebox{1.01}{\includegraphics{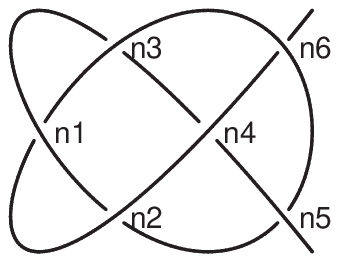}}}\\
$4_1$ &  $5_2$\\
$C(1,0,1,2)$&
%$[1,1,1,1]$&
$C(-1,-2,1,2)$
\end{tabular}
\end{center}
\caption{Knots with the Chebyshev diagram $\cC(4,5)$}\label{diag4}
\end{figure}

Once again, the situation is best explained by typical examples.
Figure \ref{diag4} represents two knots with the same  Chebyshev diagram
$\cC(4,5): x= T_4(t), \  y= T_5(t)$.  A Schubert fraction of the first knot
is $ \Frac 52 = [1,0,1,2]$, it is the figure-eight knot.
A Schubert fraction of the second knot is
$ \Frac 7 {-4} = [-1,-2,1,2]$, it is the twist knot $5_2$.
\subsection{Continued fractions}
Let $\alpha,\beta$ be relatively prime integers.
Then $ \Frac {\alpha }{\beta}$ admits the continued fraction expansion
$\Frac{\alpha}{\beta} = [q_1, q_2, \ldots, q_n ] $
if and only if there exist integers  $r_i$ such that
$$
\left \{
\begin{array}{rcl}
\alpha &=& q_1 \beta +r_2,\\
\beta &=& q_2 r_2 + r_3,\\
&\vdots&\\
r_{n-2} &=& q_{n-1} r_{n-1} + r_n,\\
r_{n-1} &=& q_n r_n.
\end{array}
\right .
$$
The integers $q_i$ are called the quotients of the continued fraction.
Euclidean algorithms provide various continued
fraction expansions which are useful to the study of two-bridge knots
(see \cite{BZ,St}).
\begin{definition}
Let $r >0$ be a rational number, and $r=[q_1, \ldots, q_n ] $ be a continued fraction
with $q_i >0.$
The {\em crossing number } of $r$ is defined by
$\cn(r)= q_1 + \cdots + q_n. $
\end{definition}
\begin{remark}
When $q_i$ are positive integers, the continued fraction
expansion $[q_1, q_2, \ldots, q_n ] $ is unique up to $[q_n]=[q_{n}-1,1]$.
$\cn(\Frac\alpha\beta)$ is the crossing number of the
knot $K= S\bigl ( \Frac\alpha\beta  \bigr) $. It means that it is  the
minimum number of crossing points for all diagrams of $K$
%When
%the quotient $q_i$ are positive, the diagram is alternate and has a minimal
%number of crossing points
(\cite{Mu}).
\end{remark}
\pn
We shall be interested by algorithms
where the sequence of remainders is not necessarily decreasing anymore.
In this case, if $\Frac\alpha\beta = [a_1, \ldots, a_n]$, we have
$\cn(\Frac\alpha\beta) \leq \sum_{k=1}^n \abs{a_i}$.
\begin{definition}
A continued fraction $ [ a_1,a_2, \ldots, a_n]$ is regular if it
has the following\\properties:
$$
a_i \neq 0, \,  a_{n-1} a_n >0,   \hbox{ and }
a_i a_{i+1} <0 \Rightarrow  a_{i+1} a_{i+2} >0, \ i = 1, \ldots, n-2.$$
If $a_1a_2>0$ she shall say that the continued fraction is biregular.
\end{definition}
\begin{proposition}\label{bireg}
Let $\Frac \alpha\beta = [a_1, \ldots, a_n]$
be a biregular continued fraction. Its crossing number is
\[
\cn(\Frac \alpha\beta) = \sum_{k=1}^n \abs{a_i} - \sharp \{i,
a_ia_{i+1}<0\}.
\label{biregf}
\]
\end{proposition}
\Pf We prove
this result by induction on the number of sign changes
$k = \sharp \{i, a_ia_{i+1}<0\}$.
If $k$ is 0, then $K$ is alternate and the result is true. If
$k>0$ let us consider the first change of sign. The Conway normal form of $K$
is $[x,a,b,-c,-d,-y]$ where $a,b,c,d$ are positive integers and $x$
is a sequence (possibly empty) of positive integers and $y$ is a sequence
of integers.
We have $ [x,a,b,-c,-d,-y] = [x,a,b-1,1,c-1,d,y]$.
\bi
\item Suppose $(b-1)(c-1)>0$, then
the sum of absolute values has decreased by 1 and the number of changes of sign
has also decreased by 1.
\item Suppose $b=1, c \not = 1$ (resp. $c=1, b\not = 1$). Then
$[x,a,b,-c,-d,-y] = [x,a,0,1,c-1,d,y] =[x,a+1,c-1,d,y]$. (resp. $[x,a-1,c+1,d,y]$).
The sum of absolute values has decreased by 1 and the number of changes of sign
has also decreased by 1.
\item Suppose $b=c=1$. Then
$[x,a,b,-c,-d,-y] = [x,a,0,1,0,d,y] =[x,a+d+1,y]$.
The sum of absolute values has decreased by 1 and the number of changes of sign
has also decreased by 1.
\EPf
\ei
Note that Formula (\ref{biregf}) is also true  in other
cases. For instance, Formula (\ref{biregf}) still holds
when $a_1, \ldots,a_n$ are non zero even integers (see \cite{St}).
\pn
We shall now use the basic (subtractive) Euclidean algorithm to get continued fractions
of the form $ [ \pm 1, \pm 1, \ldots, \pm 1 ]$.
\subsection{Continued fractions ${\mathbf{[\pm 1,\pm 1,\ldots,\pm 1]}}$}
We will consider the following homographies:
\[
P : x \mapsto [1,x] = 1+\Frac 1x, \,
M : x \mapsto [1,-1,-x] = \Frac{1}{1+x}. \label{pm}
\]
Let $E$ be the set of positive real numbers.
We have  $P(E) = ]1, \infty [$ and $M(E)= ]0, 1 [.$
%Consequently, $ P(E) \subset E, \  M(E) \subset E,$ so
$P(E)$ and $M(E)$ are disjoint subsets of $E.$
\begin{theorem}\label{th1}
Let $ \Frac \alpha\beta >0 $ be a rational number.
There is a unique regular continued fraction such that
$$\Frac \alpha\beta = [1, e_2, \ldots,  e_n], \,  e_i= \pm 1.$$
Furthermore, $\alpha>\beta$ if and only if $[e_1, e_2, \ldots, e_n]$ is biregular.
\end{theorem}
\Pf
Let us prove the existence by induction on the height
$h(\Frac\alpha\beta)=\alpha+\beta$.
\bi
\item
If $h=2$ then  $\Frac \alpha \beta = 1 = [1]$ and the  result is true.
\item
If $ \alpha > \beta,$ we have
$\Frac \alpha \beta = P(\Frac \beta {\alpha - \beta}) =
[1,\Frac \beta {\alpha - \beta }].
$
Since $h( \Frac \beta {\alpha - \beta } ) < h ( \Frac \alpha \beta )$,
we get our regular continued fraction for $r$ by induction.
\item
If $ \beta > \alpha$  we have
$ \Frac \alpha \beta = M( \Frac  {\beta - \alpha } \alpha ) =
[1, -1, -\Frac {\beta - \alpha} \alpha  ].
$
And we also get a regular continued fraction for $r.$
\ei
Conversely, let $r$ be defined by the regular continued fraction
$ r = [ 1, r_2, \ldots, r_n]$, $r_i = \pm 1$,  $n \ge 2 $.
Let us prove, by induction on the length $n$ of the continued fraction, that $r>0$
and that $r >1$ if and only if $r_2 = 1$.
\bi
\item If $r_2=1$ we have $r=P([1,r_3,\ldots,r_n] ),$ and by induction $r\in P(E)$ and
then $r>1$.
\item If $r_2=-1,$ we have $r_3=-1$ and $r= M([1, -r_4,\ldots, -r_n])$.
By induction, $r\in M(E) = ]0,1[$.
\ei
The uniqueness is now easy to prove.
Let $ r = [1, r_2, \ldots, r_n] = [1, r'_2, \ldots, r'_{n'} ].$
\bi
\item If $r>1$ then $r_2=r'_2=1$ and
$[1,1, r_3, \ldots, r_n ] =[1, 1, r'_3, \ldots r'_{n'} ]$.
Consequently,
$[1, r_3, \ldots, r_n] = [1, r'_3, \ldots, r'_{n'} ],$ and by induction
$r_i=r'_i$ for all $i.$
\item If $r < 1,$ then $r_2=r_3=r'_2=r'_3=-1$ and
$[1,-1,-1, r_4, \ldots, r_n] =[1,-1,-1, r'_4,\ldots, r'_{n'} ].$
Then,
$ [ 1, -r_4, \ldots, -r_n] = [ 1, -r'_4, \ldots, -r'_{n'} ]$
and by induction
$ r_i=r'_i$ for all $i.$
\EPf
\ei
\begin{definition}
Let $\Frac\alpha\beta>0$ be the regular continued fraction $[e_1, \ldots, e_n]$, $e_i=\pm 1$.
We will denote its length $n$ by $\ell(\Frac\alpha\beta)$. Note that
$\ell(-\Frac\alpha\beta) = \ell(\Frac\alpha\beta)$.
\end{definition}
\begin{examples}\label{9297}
Using our algorithm we obtain
$$
\begin{array}{rcl}
\Frac 97&=& [1,\Frac 72] = [1,1,\Frac 25] =  [1,1,1,-1, -\Frac 32] =
[1,1,1,-1, -1, -\Frac 21] = [1,1,1,-1, -1, -1, -1],\\[10pt]
\Frac 92&=& [1,\Frac 27] = [1,1,-1,-\Frac 52] = [1,1,-1,-1,-\Frac 23] =
[1,1,-1,-1,-1,1,\Frac 12] \\
&=& [1,1,-1,-1,-1,1,1,-1,-1] = [4,2].
\end{array}
$$
\begin{figure}[th]
\begin{center}
\begin{tabular}{ccc}
{\scalebox{.6}{\includegraphics{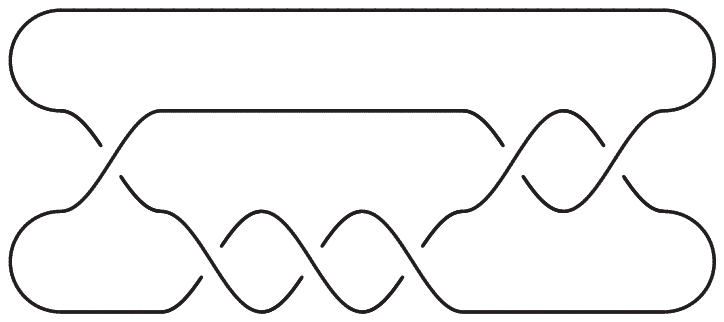}}}&
{\scalebox{.6}{\includegraphics{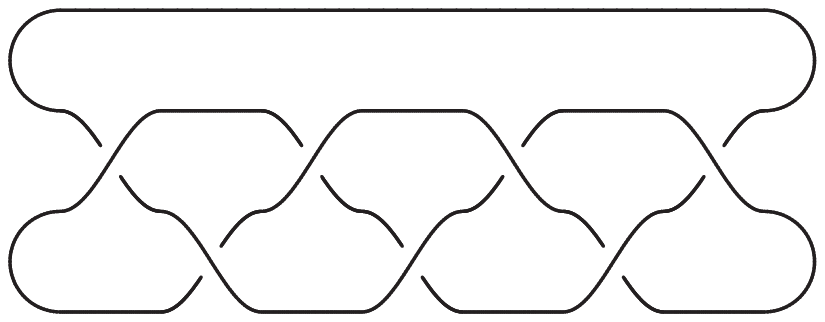}}}&
{\scalebox{.6}{\includegraphics{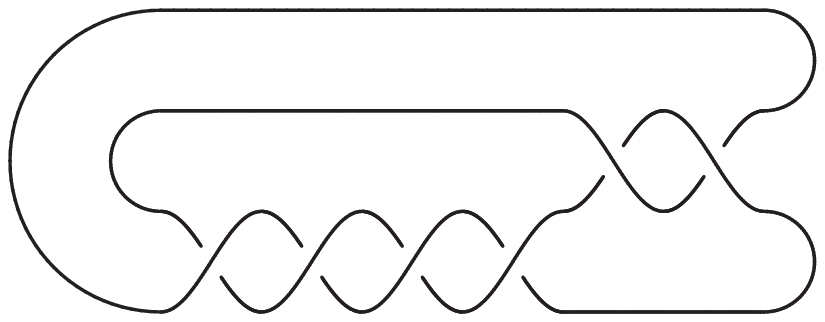}}}\\
$C(1,3,2)$&$C(1,1,1,-1,-1,-1,-1)$&$C(4,2)$
\end{tabular}
\end{center}
\caption{Diagrams of the knot $6_1 = S(\Frac 97)$ and its mirror image $S(\Frac 92)$ }
\label{k61}
\end{figure}
\end{examples}
We will rather use the notation
$$
\Frac 97 = P^2 M P^3 (\infty), \,
\Frac 92 = P MP M^2 P (\infty).
$$
We get $\ell (\Frac 97) = 7$, $\ell (\Frac 92) = 9$.
The crossing numbers of these fractions are
$\cn(\Frac 97) = \cn([1,3,2]) = 6 = 7-1$ and $\cn(\Frac 92) = \cn([4,2]) = 6 = 9-3$.
If the fractions $\Frac 97$ and $\Frac 92$ have the same crossing number,
it is because the knot $S(\Frac 97)$ is the mirror image of $S(\Frac 92)$.
In order to get a full description of two-bridge knots we shall need a more detailed
study of the homographies $P$ and $M$.
% defined by Formula \ref{pm}.
\pn
\begin{proposition}
The multiplicative monoid $\GG=\langle P,M \rangle$
is free.
%There is a bijection $g$ from the set of positive rational numbers
%to $\GG$ such that $r = g(r)(1)$.
The mapping $g: G \mapsto G(\infty)$ is a bijection from $\GG \cdot P$ to $\QQ_{>0}$ and
$g(P\cdot \GG \cdot P) = \QQ_{>1}$, the set of rational  numbers greater than $1.$
\end{proposition}
\Pf Suppose  that $PX=MX'$ for some $X,X'$ in $\GG$. Then we
would have $PX(1)=MX'(1) \in P(E)\bigcap M(E) = \emptyset$. Clearly, this
means that $\GG$ is free. Similarly, from $P(\infty)=1$, we deduce that the
mapping $G \mapsto G \cdot P (\infty)$ is injective. From
Theorem \ref{th1} and $P(\infty)=1$,  we deduce that $g$ is
surjective. \EPf
\begin{remark}\label{rl}
Let $r=G(\infty)=[e_1, \ldots, e_n], \  e_i=\pm1, $ be a regular continued fraction.
It is easy to find the unique homography $G \in \GG\cdot P$ such
that $r=G(\infty)$. Consider the sequence $(e_1, \ldots, e_n)$. For any
$i$ such that $e_ie_{i+1}<0$, replace the couple $(e_i,e_{i+1})$ by
$M,$ and then replace each  remaining $e_i$ by $P$.
\pn
Let $G= P^{p_1} M^{m_1} \cdots M^{m_k} P^{p_{k+1}}$.
Let $p= p_1 + \cdots + p_{k+1}$ be the degree
of $G$ in $P$ and $m = m_1 + \cdots + m_k$ its degree in $M$. Then we have
$n = \ell(r) = p + 2m$ and $\cn(r) = p+m$.
\end{remark}
We shall consider matrix notations for many proofs. We will consider
$$\cmatrix{\alpha\cr\beta} = P^{p_1} M^{m_1} \cdots M^{m_k} P^{p_{k+1}}
\cmatrix{1\cr 0}, \quad
P = \cmatrix{1&1\cr1&0},\quad M =\cmatrix{0&1\cr1&1}.
$$
\begin{lemma}\label{length}
Let $\Frac{\alpha}{\beta} = [e_1, \ldots, e_n]$ be a regular continued
fraction ($e_i = \pm 1$). We have
\bi
\item $n\equiv 2 \Mod 3$ if and only if $\alpha$ is even and $\beta$ is odd.
\item $n\equiv 0 \Mod 3$ if and only if $\alpha$ is odd and $\beta$ is even.
\item $n\equiv 1 \Mod 3$ if and only if $\alpha$ and $\beta$ are odd.
\ei
\end{lemma}
\Pf
Let us write  $\cmatrix{\alpha\cr\beta} = P^{p_1} M^{m_1} \cdots M^{m_k} P^{p_{k+1}}
\cmatrix{1\cr 0}$.
We have (from remark \ref{rl}) $n= p+2m.$
Since $P^3 \equiv \id $ and $M \equiv P^2 \Mod 2$, we get
$P^{p_1} M^{m_1} \cdots M^{m_k} P^{p_{k+1}} \equiv P^n \Mod 2$
\bi
\item[]
If $n\equiv 2 \Mod 3$, then $\cmatrix{\alpha\cr \beta} \equiv M \cmatrix{1\cr 0} \equiv \cmatrix{0\cr 1} \Mod 2$.
\item[]
If $n\equiv 1 \Mod 3$, then $\cmatrix{\alpha\cr \beta} \equiv P \cmatrix{1\cr 0} \equiv \cmatrix{1\cr 1} \Mod 2$.
\item[]
If $n\equiv 0 \Mod 3$, then $\cmatrix{\alpha\cr \beta} \equiv \cmatrix{1\cr 0} \Mod 2$.
\EPf
\ei
%%%%%%%%%%%%%%%%%%%%%%%%%%%%%%%%%%%%%%%%%%%%%%%%%%%
%In terms of vectors in $\ZZ^2$, we get for example
%$\cmatrix{9\cr 7} = P^2 M P^3 \cmatrix{1\cr 0}, \
%\cmatrix{9\cr 2} = P M P M^2 P \cmatrix{1\cr 0}$.
%\pn
\begin{definition}\\
We define on $\GG$ the anti-homomorphism $G\mapsto \overline G$ with $\overline M=M$, $\overline P=P$. \\
We define on $\GG$ the homomorphism $G \mapsto \hat G$ with $\hat M = P, \, \hat P = M$.
\end{definition}
\begin{proposition}\label{gg}
Let $\alpha>\beta$ and consider $\Frac\alpha\beta=PGP(\infty)$ and $N=\cn(\Frac\alpha\beta)$.
Let $\beta'$ such that $0<\beta'<\alpha$ and $\beta\beta'\equiv (-1)^{N-1}\Mod\alpha$.
Then we have
$$
\Frac\beta\alpha = M\hat G P(\infty),\quad
\Frac\alpha{\alpha-\beta} = P \hat G P(\infty), \quad
\Frac\alpha{\beta'} = P \overline G P(\infty).$$
We also have
$$
\ell (\Frac {\beta}{\alpha}) + \ell( \Frac {\alpha}{\beta }) = 3 N -1, \quad
\ell ( \Frac {\alpha}{\alpha - \beta } ) + \ell ( \Frac {\alpha}{\beta})
= 3 N -2,
\quad
\ell ( \Frac {\alpha}{\beta'} ) = \ell( \Frac {\alpha}{\beta }).
$$
\end{proposition}
\Pf
We use matrix notations for this proof.
%\pn
Let us consider
$$PGP=\cmatrix{\alpha&\beta'\cr\beta&\alpha'}= P^{p_1} M^{m_1} \cdots M^{m_k} P^{p_{k+1}}.$$
From $\det P = \det M  = -1$, we obtain
$\alpha\alpha' - \beta\beta' = (-1)^N$.
Let $A =\cmatrix{a&c\cr b&d}$ be a matrix such that
$0\leq c \leq a$ and $0\leq d \leq b$.
From $PA = \cmatrix{a+c&b+d\cr a&c}$ and
$MA = \cmatrix{b&d\cr a+b&c+d}$ we deduce that
$PGP$ satisfies $0<\alpha'<\beta$ and $0<\beta'<\alpha$.
We therefore conclude that,
$\beta'$ is the integer defined by
$0<\beta'<\alpha$, $\beta\beta'\equiv (-1)^{N-1} \Mod\alpha$.
By transposition we deduce that
$$\cmatrix{\alpha&\beta\cr\beta'&\alpha'} =
P^{p_{k+1}} M^{m_{k}} \cdots M^{m_1} P^{p_1} = P \overline G P,$$
which implies $\Frac{\alpha}{\beta'} = P \overline G P(\infty)$.
\pn
Let us introduce $J = \cmatrix{0&1\cr1&0}$. We have $J^2 = \id, M=JPJ$ and $P=JMJ.$
Therefore
$$
\cmatrix{\beta\cr\alpha} = J \cmatrix{\alpha\cr\beta}=
M^{p_1} P^{m_1} \cdots P^{m_k} M^{p_{k+1}-1} J P \cmatrix{1 \cr 0}=
M^{p_1} P^{m_1} \cdots P^{m_k} M^{p_{k+1}-1} P \cmatrix{1 \cr 0},
$$
that is $\Frac\beta\alpha = M\hat G P (\infty) $.
\pn
Then,
$\cmatrix{\alpha \cr \alpha - \beta } = PM^{-1} \cmatrix{ \beta \cr \alpha }
= P \hat G P \cmatrix{1 \cr 0}.$
That is
$ \Frac {\alpha}{ \alpha - \beta} = P \hat G P ( \infty ).$
\pn
Relations on lengths are derived from the previous relations and remark \ref{rl}.
\EPf
%%%%%%%%%%%%%%%%%%%%%%%%%%%%%%%%%%%%%%%%%%%%%%%%%%%%%%%%%%%%%%%%%%%%%%%%%%%%
\pn
We deduce
\begin{proposition}\label{palin}
Let $G \in \GG $ and  $\Frac\alpha\beta = [e_1, \ldots, e_n] = PGP(\infty).$
Let $K = S( \Frac \alpha\beta )$ and $N= \cn (K)$.
The following properties are equivalent:
\bn
\item $G$ is palindromic (i.e. $\overline G = G$).
\item the sequence of sign changes in $[e_1, \ldots, e_n]$ is palindromic
(i.e. $e_ie_{i+1} = e_{n-i}e_{n-i+1}$).
\item $\beta^2 \equiv (-1)^{N-1}\Mod \alpha$.
\en
Moreover we have
\bi
\item
$\beta^2 \equiv -1 \Mod \alpha$ (i.e. $K= \overline{K}$ is amphicheiral)
if and only if $N$ is even and $G = \overline G$. Furthermore, the length
$n=\ell( \frac \alpha \beta) $ is even and the sequence $[e_1, \ldots, e_n]$
is palindromic (i.e. $e_i = e_{n-i+1}$).
\item
$\beta^2\equiv 1 \Mod\alpha$ if and only if $N$ is odd and $G = \overline G$
or $N$ is even and $\hat G = {\overline G}$ (in this case $K$ is a two-component link).
\ei
\end{proposition}
\Pf
From Remark \ref{rl}, it is straightforward that if $[\eps_1, \ldots, \eps_n] = PGP(\infty)$ then
$P\overline G P(\infty) = \eps_n [\eps_n, \ldots, \eps_1]$. We deduce that
$G=\overline G$ is palindromic if and and only if the sequence of sign changes
in $[e_1, \ldots, e_n]$ is palindromic.

Let $0<\beta'<\alpha$ such that $\beta'\beta \equiv (-1)^{N-1}\Mod\alpha$.
We have from the previous proposition: $\Frac\alpha{\beta'} = P \overline{G}P ( \infty )$
We thus deduce that  $G=\overline G$ is equivalent to $\beta = \beta'$, that is
$\beta^2\equiv(-1)^{N-1}\Mod\alpha$.

Suppose now that $\beta^2\equiv 1\Mod\alpha$. If $N$ is even then $\beta'=\alpha-\beta$,
that is $P\overline GP = P\hat G P$ and $\overline G=\hat G$.
We  have $p+2m = m + 2 p-2$ and then $2n =2(p+2m) = 3N-2$. This implies
$n\equiv 2\Mod 3$. By  Lemma \ref{length}, $\alpha$ is even and $K$ is a two-component link.
If $N$ is odd then $\beta'=\beta$ and $G=\overline G$ by the first part of our proof.

Suppose now that $\beta^2\equiv -1 \Mod\alpha$.
If $N$ is odd then $\beta'=\alpha-\beta$ and
by the same argument we should have $n=3N-n-2,$
which would imply that $N$ is even. We deduce
that amphicheiral rational links have even crossing numbers
and from $\beta'=\beta$ we get $G=\overline G$.
The crossing number $N=m+p$ is even and $G$ is palindromic so $m$ and $p$ are both even.
Consequently $n= p+2m$ is even and the number of sign changes is even. We thus have
$e_n=1$ and $(e_n, \ldots, e_1)=(e_1, \ldots, e_n)$.
\EPf
\pn
We deduce also a method to find a minimal Chebyshev diagram for any rational knot.
\begin{proposition}\label{ell}
Let $K$ be a two-bridge knot with crossing number $N$.
\bn
\item There exists $\Frac{\alpha}{\beta}>1$ such that
$K=S(\pm \Frac\alpha\beta)$ and $n=\ell(\Frac{\alpha}{\beta}) < \frac 32 N-1$.
\item There exists a biregular sequence
$(e_1, \ldots, e_n)$, $e_i=\pm 1$, such that
$K=C(e_1,\ldots,e_n)$.% is a Conway normal form for $K$.
%\item There exists $b<\Frac 32 N$, such that $x=T_3(t),\, y=T_b(t)$ is a Chebyshev diagram for $K$.
\item If $K=C(\eps_1,\ldots,\eps_m)$, $\eps_i=\pm1$, then $m\geq n\geq N$.
\en
\end{proposition}
Let $K$ be a two-bridge knot.
Let $r=\Frac\alpha\beta>1$ such that $K=S(r)$.
We have $\overline K = S(r')$ where $r'=\Frac{\alpha}{\alpha-\beta}$.
From Proposition \ref{gg} and Proposition \ref{bireg}, we have
$\ell(r)+\ell(r')=3N-2$ and therefore $N\leq \min(\ell(r),\ell(r')) < \frac 32 N$.
From Lemma \ref{length}, we have $\ell(r)\not \equiv 2 \Mod 3$ so $\ell(r)\not=\ell(r')$
and $\min(\ell(r),\ell(r'))<\frac 32 N -1$.
We may suppose that $n=\ell(r)$. Let
$r= \pm [e_1, \ldots, e_n]$, then $C(e_1,\ldots,e_n)$ is a Conway normal
form for $K$ with $n<\frac 32 N -1$. This Conway normal form is a Chebyshev diagram
$\cC(3,n+1): x= T_3(t),\, y=T_{n+1}(t)$.
\pn
Let us consider $\gamma$ such that
$\beta\gamma \equiv 1 \Mod\alpha$ and $0<\gamma<\alpha$.
Let $\rho = \Frac{\alpha}{\gamma}$ and $\rho' = \Frac{\alpha}{\alpha-\gamma}$.
We have $K = S(\rho)$ and $\overline K = S(\rho')$ and from Proposition \ref{gg}:
$\min(\ell(\rho),\ell(\rho'))=\min(\ell(r),\ell(r'))$.
Suppose that $K=C(\eps_1, \ldots, \eps_m)$ then
$K = S(x)$ where $x=[\eps_1, \ldots, \eps_m]$.
We thus deduce that $x=\Frac{\alpha}{\beta+k\alpha}$ or $x' = \Frac{\alpha}{\gamma+k\alpha}$ where $k \in \ZZ$.
\bn
\item[--] If $k=2p>0$ then
we have $x=(MP)^p r$ so $m=\ell(x)=\ell(r)+3p>\ell(r)$.
\item[--] If $k=2p+1>0$ then
$x=(MP)^p M (\Frac 1r )$ so $\ell(x)=\ell(1/r)+3p+2=\ell(r')+3p+3>\ell(r')$.
\item[--] If $k=-(2p+1)<0$ then $-x=(MP)^p (r')$ so $\ell(x)=\ell(-x)=\ell(r')+3p>\ell(r')$.
\item[--] If $k=-2p>0$ then
$-x=(MP)^{p-1} M (\Frac 1{r'})$ so $\ell(x)=\ell(1/r')+3p-1=\ell(r)+3p>\ell(r)$.
\en
If $x'=\Frac{\alpha}{\gamma+k\alpha}$, we obtain the same relations.
We deduce that $m\geq\min(\ell(r),\ell(r'))$.
\EPf
%%%%%%%%%%%%%%%%%%%%%%%%%%%%%%%%%%%%%%%%%%%%%%%%%%%%%%%%%%%%%%%%
\begin{remark}[Computing the minimal Chebyshev diagram $\mathbf{\cC(3,b)}$]\label{minb}\\
Let $K=S(\Frac\alpha\beta)$ with $\Frac\alpha\beta = PGP(\infty)$.
The condition $\ell(\Frac\alpha\beta)$ is minimal, that is
$\ell(\Frac\alpha\beta) < \frac 32 N -1$, is equivalent to
$p\geq m+3$ where $p=\deg_P(PGP)$ and $m=\deg_M(PGP)$.
In this case $b=\ell(\Frac\alpha\beta)+1$ is the smallest integer such
that $K=S(\Frac\alpha\beta)$ has a Chebyshev diagram $x=T_3(t),\, y=T_b(t)$.
If $p<m+3$, using Proposition \ref{gg}, $K$ has a Chebyshev diagram
$\cC(3,b')$ with $b' = 3N - b < \frac 32 N$. This last diagram is minimal.
\end{remark}
\begin{example}[Torus knots]\label{tok}
The Schubert fraction of the torus knot $T(2,2k+1)$ is $2k+1$.
We have $ PM(x) = x+2,$ and then $ (PM)^k (x) = x+2k,  (PM)^kP (x)= 2k+1+ \Frac 1x.$
So we get the continued fraction of length $3k+1$:
$ 2k+1 = (PM)^kP(\infty)$.
This shows that the torus knot $T(2, 2k+1)$
has a Chebyshev diagram $ \cC( 3, 3k+2).$
This is not a minimal diagram.
\pn
On the other hand, we get $(PM)^{k-1}P^2(\infty) = 2k$ so
$\Frac{2k+1}{2k} = P(PM)^{k-1}P^2(\infty) >1$.
This shows that the torus knot $T(2, 2k+1)$
has a Chebyshev diagram $ \cC( 3,3k+1).$
This diagram is minimal by Remark \ref{minb}.
We will see that $T(2, 2k+1)$  is in fact a harmonic knot.
\end{example}
\begin{example}[Twist knots]\label{twk}
The twist knot $\cT_{n}$ is defined by $\cT_n = S(n+\Frac{1}{2})$.
%Its crossing number is therefore $n+2$.
\pn
From $P^3( x) = \Frac{3x+2}{2x+1}$, we get the continued fraction of length $3k+3$:
$\Frac {4k+3}2= (PM)^kP^3 (\infty).$
This shows that the twist knot $\cT_{2k+1}$
has a  Chebyshev diagram $\cC ( 3, 3k+4),$ which is minimal by Remark \ref{minb}.
\pn
We also deduce that $P(PM)^{k-1}P^3 (\infty) = P \left(\frac 12 (4k-1)\right ) = \Frac{4k+1}{4k-1}$.
This shows that the twist knot $\cT_{2k}$ has a  minimal Chebyshev diagram $\cC (3,{3k+2}).$
\pn
We shall see that these knots are not harmonic knots for $a=3$
and we will give explicit bounds for their polynomial parametrizations.
\end{example}
\begin{example}[Generalized stevedore knots]\label{sk3}
The generalized stevedore knot $\cS_k$ is defined by $\cS_k = S(2k+2+\Frac{1}{2k})$.
We have
$$
(MP)^k = \cmatrix{1&0\cr 2k&1}, \, (PM)^k =\cmatrix{1&2k\cr 0&1}
$$
so $2k+2+\Frac{1}{2k} = (PM)^{k+1}(MP)^k (\infty)$. This shows
that the stevedore knots have a Chebyshev diagram $\cC(3,6k+4)$.
It is not minimal and
we see, using Remark \ref{minb}, that the
stevedore knot $\cS_k$  also has a minimal Chebyshev diagram $\cC(3,6k+2)$.
Moreover, using Proposition \ref{gg}, we get
$$
\Frac{(k+1)^2}{(k+1)^2-2k} = P^2 (MP)^k (PM)^{k-1} P^2 (\infty).
$$
\end{example}
%\subsection{Enumeration of rational links}
\subsection{Continued fractions ${\mathbf{[\pm 1,\pm 2,\ldots,\pm 1,\pm 2]}}$}
%In the study of harmonic knots of the form $H(4,b,c),$ we shall need the following result
Let us consider the homographies $A(x)=[1,2,x] =\Frac{3x+1}{2x+1}$,
$B(x)=[1,-2,-x]=\Frac{x+1}{2x+1}$, $S(x)=-x$.
We shall also use the classical matrix notation for these homographies
$$
A = \cmatrix{3&1\cr2&1}, \,
B = \cmatrix{1&1\cr2&1}, \,
S = \cmatrix{1&0\cr0&-1}.
$$
\begin{lemma}\label{nn}
The monoid $\Ga = \langle (AS)^k A, (AS)^{k+1}B,
B(SA)^k,   \  B(SA)^{k+1}SB, \  k\geq 0\rangle$
%The elements of $\Ga$ have nonnegative coefficients.
%The only element with $G\infty = \infty$ is $ASB$.
%The monoid $\Ga $
is free.
Let $ \Ga ^* $ be the subset of elements of $ \Ga$ that are not of the form
$ M \cdot (AS)^{k+1} B, $ or $ M\cdot B(SA)^{k+1}SB$.
There is an injection $h : \Ga ^* \to \QQ_{>0}$ such that
$h(G) = G(\infty)$.
\end{lemma}
\Pf
Let us denote $E= \RR ^* _+ \bigcup \{ \infty  \}= ]0, \infty].$
We will describe $G(E)$ for any generator $G$ of $ \Ga.$
Let $C=AS =\cmatrix{3&-1\cr2&-1}$. We get $C^2 = 2C+\id,$ and then
%$C^{k+2} = 2C^{k+1}+C^k$. Then we have
$C^{k+2}A = 2C^{k+1}A+C^kA$.
From $CA = ASA =\cmatrix{7&2\cr4&1}$ we deduce by induction that
%$C^k A$ has positive coefficients and
$ (AS)^k A (E) \subset ]1, \infty [ .$
Similarly, we get
$ (AS)^{k+1} B(E) \subset ]1, \infty], \  B(SA)^k(E)\subset [ 0, 1[, $ and
$ B( SA) ^{k+1}SB (E) \subset ]0, 1[ .$
%We also remark that
%$ASB=CB= \cmatrix{ 1 & 2 \cr 0 & 1 }$ is the only generator of $\Ga$ such
%that $G\cmatrix{1\cr0} = \cmatrix{1\cr0}$.
\pn
Now, we shall prove that if $G$ and $G'$ are distinct generators of $ \Ga,$
the relation $ G \cdot M = G' \cdot M' $ is impossible.
This is clear in cases where  $G(E)$ and $G'(E)$ are disjoint.
Suppose that  $ ( AS)^k AM= (AS)^{k'} BM' .$
If $k<k',$ we would get
$ M= S \Bigl( ( AS) ^{k''} BM' \Bigr) = S M'' ,$ which is impossible
because $ S M''$ has some negative coefficients.
If $ k \ge k'$, we get $ (AS)^{k''} AM= BM'$, which is also impossible because
$ B(E)$ and $ ( AS)^{k''} A (E) $ are disjoint.
The proof of the impossibility of
$ B (SA)^k M= B(SA)^{k'+1} SB M' $ is analogous to the preceding one. This completes
the proof that $ \Ga$ is free.
\pn
Now, let us prove that $ M \neq M', \  M, M' \in \Ga  ^*$ implies $ M(\infty) \neq M'(\infty).$
There are two cases that are not obvious.

If $ (AS)^kA M(\infty)= (AS)^{k'} BM' ( \infty) , \   M, M' \in \Ga ^* .$
If $ k< k'$, we get $ M(\infty) = S ( M'' ( \infty)) ,$ which is impossible since
$ M(\infty) \in ]0, \infty[, $ and $ SM'' (\infty) \in ] -\infty , 0[ .$\\
If $ k \ge k',$ we get
$ BM'(\infty) = ( AS)^{k''} A M (\infty) ,$
which is also impossible because
$BM'(\infty) <1$ and $ (AS)^{k''} AM( \infty) >1.$

The proof of the impossibility of
$B (SA)^k M (\infty) = B(SA)^{k'+1} SB M' (\infty) $ is analogous to the preceding one.
This completes the proof of the injectivity of $ h.$
\EPf
\begin{theorem}  \label{cf1212}
Let $r=\Frac\alpha\beta>0$ be a rational number.
Then $r$ has a continued fraction expansion
$\Frac\alpha\beta = [1,2 e_2,\ldots,e_{2n-1},2 e_{2n}]$, $e_i = \pm 1$,
if and only if $\beta$ is even.
Furthermore, we can suppose that there are no three consecutive sign changes.
In this case the continued fraction is unique, and $\alpha > \beta$ if and only
if $e_1= e_2= 1$.
\end{theorem}
\Pf
Let us suppose that $\Frac\alpha\beta =[1,\pm2,\ldots,\pm1,\pm2]$.
It means that $\cmatrix{\alpha\cr\beta} = G \cmatrix{1\cr0}$
where $G \in \langle A,B,S\rangle$.
But we have
$A\equiv B \equiv \cmatrix{1&1\cr0&1} \Mod 2$, $S \equiv \id \Mod 2$.
Consequently,
$\cmatrix{\alpha\cr\beta}
\equiv \cmatrix{1&1\cr0&1}^k \cmatrix{1\cr0}
\equiv \cmatrix{1\cr0}\Mod 2$ and $\beta$ is even.
\pn
Let us suppose that $\beta$ is even.
We shall use induction on the height
$h(\Frac\alpha\beta) = \max(\abs\alpha,\abs\beta)$.
\bi
\item If $h(\Frac\alpha\beta)=2,$
then $\alpha=1$ and $\beta=2$ and we have $r=[1,-2]=B(\infty)$.
\item We have two cases to consider
\bn
\item[--] If $\alpha>\beta>0$ then we write
$
\Frac\alpha\beta = [1,2,-1,2,\Frac{\alpha-2\beta}{\beta}]=
ASB \left (\Frac{\alpha-2\beta}{\beta}\right ).
$
We have $h(\Frac{\alpha-2\beta}{\beta}) < h(\Frac\alpha\beta)$
and we conclude by induction.
\item[--] If $\beta > \alpha >0$ we write
$
\Frac\alpha\beta = [1,-2,\Frac{\alpha-\beta}{2\alpha-\beta}]=
B \left (\Frac{\beta - \alpha}{2\alpha-\beta}\right ).
$
From $\abs{2\beta-\alpha}\leq \beta$ we have
$h(\Frac{\alpha-\beta}{2\alpha-\beta})<h(\Frac\alpha\beta)$ and we conclude
by induction.
\en
\ei
The existence of a continued fraction
$[ 1, \pm 2, \ldots, \pm 1,  \pm2]$ is proved.
\pn
Using the identities $BSBS(x)=[1,-2,1,-2,x]=x$ and $[2,-1,2,-1,x]=x,$
we can delete all subsequences
with three consecutive sign changes in our sequence.
We deduce also that
$\cmatrix{\alpha\cr\beta} = G \cmatrix{1\cr 0}$ where
$G \in \langle A,S,B\rangle$.
Furthermore, from $BSB(x) = [1,-2,1,-2,-x]$
we see that $G$ contains no $BSB$.
We also see from $ASBSA(x) = [1,2,-1,2,-1,-2,-x]$ that
$G$ contains no $SBS$.
Consequently $G $ is an element of $ \Ga ^* $
%is a product of matrices
%$(AS)^kA, \, (AS)^k B, \, B(SA)^k, \, k\geq 0$
%only
and, by Lemma \ref{nn},  the continued fraction
$\Frac\alpha\beta = [1,\pm2,\ldots,\pm1,\pm2]$ is unique.
\EPf
%%%%%%%%%%%%%%%%%%%%%%%%%%%%%%%%%%%%%%%%%%%%%%%%%%%%%%%%%%
\begin{example}[Torus knots]\label{tnh4}\\
Since $BSA = \cmatrix{1&0\cr4&1}, \ $ we get by induction
$$
A (BSA)^k = \cmatrix{4k+3&1\cr 4k+2&1}, \,
A (BSA)^k B = \cmatrix{4k+5&4k+4\cr 4k+4&4k+3}.
$$
We deduce  the following continued fractions
$$ \Frac {4k+3}{4k+2}= [ 1,2, \underline{1,-2,1,2}, \ldots, \underline{1, -2,1,2} ], \   \quad
\Frac {4k+5}{4k+4}= [1,2, \underline{1, -2, 1,2}, \ldots, \underline{1,-2,1,2}, 1, -2 ] $$
of length $4k+2$ and $4k+4.$
Since the knot with Schubert fraction
$ \Frac {N}{N-1}$ is the torus knot $\overline{T}(2,N),$
we see that this knot admits a Chebyshev diagram
$ x= T_4(t), y= T_N (t).$
\end{example}
\begin{example}[Twist knots]\label{twk4}\\
The twist knot is the knot
$\cT_m= S( \Frac {2m+1} 2 ) = S ( \Frac {2m+1}{m+1} )= \overline{S}(\Frac{2m+1}m).$
We have the  continued fractions
$$ \Frac{8k+1}{4k} = (ASB) (BSA)^k ( \infty ), \ \quad
\Frac{8k+5}{4k+2} = (ASB)( BSA)^k B( \infty ).
$$
We deduce that $\cT_{2n}$ has a Chebyshev diagram
$ x= T_4(t), \, y= T_ {2n+5}(t).$
We get similarly
$$
\Frac {8k+7}{4k+4} = ASA (BSA)^k ( \infty ), \
\Frac {8k+3}{4k+2} = ASA (BSA)^{k-1} B( \infty ), \
$$
and we  deduce that   $\cT_{2n+1}$ has a Chebyshev diagram
$x= T_4(t), \,  y= T_{2n+3} (t).$
\end{example}
\begin{example}[Generalized stevedore knots]\label{sk4}\\
The generalized stevedore knot $\cS_m$ is defined by
$\cS_m= S(2m+2 + \Frac 1{2m}), \  m \ge 1.$
The stevedore knot is  $\cS_1= \overline{6}_1$.
We have $ \overline{\cS}_m = S \Bigl(\Frac{(2m+1)^2}{2m+2}\Bigr) $
and the continued fractions
$$ \Frac {(4k+1)^2}{4k+2} = (ASB)^{2k} (BSA)^k B (\infty), \
\Frac {(4k-1)^2}{4k} = (ASB)^{2k-1}(BSA)^k (\infty ).
$$
Consequently,  the stevedore knot $\cS_m$ has a Chebyshev diagram
$\cC(4,6m+3)$.% x= T_4(t), \,  y = T _{6m+3}(t).$
\end{example}
These examples show that our continued fractions are not necessarily regular.
In fact, the subsequences $ \pm (2,-1,2) $ of the continued
fraction correspond bijectively to factors $(ASB).$
\pn
There is a formula to compute the crossing number of such knots.
\begin{proposition}\label{ilot}
Let $\Frac{\alpha}{\beta} = [e_1,2e_2,\ldots,e_{2n-1},2e_{2n}]$, $e_i=\pm 1$,
where there are no three consecutive sign changes and $e_1=e_2$.
We say that $i$ is an islet in $[a_1,a_2,\ldots,a_{n}]$ when
$$\abs{a_i}=1, \, a_{i-1}=a_{i+1}=-2a_i.$$
We denote by $\sigma$ the number of islets. We have
\[
\cn(\Frac\alpha\beta) =
\sum_{k=1}^n \abs{a_i} - \sharp \{i, a_ia_{i+1}<0\} - 2 \sigma.
\label{ilotf}
\]
\end{proposition}
\Pf
By induction on the number
of double sign changes $k=\sharp\{i, e_{i-1}e_{i}<0, \, e_{i}e_{i+1}<0\}$.
If $k=0$, the sequence is biregular and
Formula (\ref{ilotf}) is true by Proposition \ref{bireg}.
Suppose $k\geq 1$.
First, we have  $a_1 a_2 >0$.
Then $\Frac\alpha\beta =\pm [x,a,b,-c,d,e,y]$ where $[x,a,b]$ is a biregular
sequence and $a,b,c,d,e>0$.
We have $[x,a,b,-c,d,e,y] = [x,a,b-1,1,c-1,-d,-e,-y] $.
\bi
\item If $c = 2$ (there is no islet at $c$), then we have $b=d=1$ and
$[x,a,b,-c,d,e,y] = [x,a+1,c-1,-d,-e,-y]$.
The sum of the absolute values has decreased by 1, as has the number of
sign changes, $\sigma$ is unchanged.
\item If $c=1$ (there is an islet at $c$) then
$[x,a,b,-c,d,e,y] = [x,a,b-1,-(d-1),-e,-y]$.
The sum of absolute values has decreased by 3, the number of sign changes by 1, $\sigma$ by 1.
\ei
In both cases,
$\sum_{k=1}^n \abs{a_i} - \sharp \{i,a_ia_{i+1}<0\} - 2\sharp\{i, e_{2i}e_{2i+1}<0, \, e_{2i+1}e_{2i+2}<0\}$
remains unchanged while $k$ has decreased by 1.
%We obtain a new sequence for which
%double sign changes occur only as $[b,-c,d]$ with $c=\pm 1$ and $b=d=\mp 2$ or
%$c=\pm 2$ and $b=d=\mp 1$.
\EPf
\pn
We shall need a specific result for biregular fractions $[\pm 1,\pm 2,\ldots, \pm 1,\pm 2]$.
\begin{proposition} \label{palin4}
Let $r= \Frac \alpha \beta$ be a rational number given by a
biregular continued fraction of the form
$ r= [e_1,  2e_2, e_3,  2e_4, \ldots e_{2m-1}, 2 e_{2m} ], \, e_1=1, \,    e_i= \pm 1.$
If the sequence of sign changes is palindromic,
i.e. if $e_k e_{k+1} = e_{2m-k} e_{2m-k+1},$
we have $\beta^2 \equiv \pm 2 \Mod \alpha.$
\end{proposition}
\Pf
From  Theorem \ref{cf1212},
and because $\Frac\alpha\beta =  [ 1,  2e_2, e_3,  2e_4, \ldots e_{2m-1}, 2 e_{2m} ]$ is
regular, we have $\Frac\alpha\beta = G (\infty)$ where
$G \in \langle B, (AS)^k A, \, k \geq 0 \rangle \subset \Ga$.
\pn
We shall consider the mapping (analogous to matrix transposition)
$$\phi: \cmatrix {a & b \cr c & d }\mapsto \cmatrix { a & \Frac c2 \cr
  2b & d }.$$
We have $ \phi(XY)= \phi (Y) \phi (X), \ \phi(A)=A,   \  \phi(B)= B$
and $\phi((AS)^k A) = (AS)^k A$.
\pn
Let us show that $G$ is a palindromic product of terms $A_k= (AS)^k A$ and $B$.

By induction on $s = \sharp\{i, e_{i}e_{i+1}<0\}.$ If $s=0$ then $G = A^m$.
Let $k = \min \{i, e_{i}=-e_{i+1}\}$. If $k=2p$ then
$G =  A^p G'$ and $G'\in S \cdot \Ga$. We have $e_1 = \ldots = e_{2p}$ and
$e_{2(m-p)+1} = \ldots = e_{2m}$ that is $G = A^p S G' S A^p$.
The subsequence
$(-e_{2p+1}, \ldots, -e_{2m-2p})$ is still palindromic and
corresponds to $G'(\infty)$.
By regularity $e_{2p+1}= e_{2p+2},$ and
we conclude by induction. If $k=2p+1$ we have $G=A^p B G' B A^p$ and we conclude also by induction.
\pn
Hence $ \phi (G)= G,$ and since $G\cmatrix {1 \cr 0 } = \cmatrix {\alpha \cr \beta },$
we see that $G$ has the form $G= \cmatrix{ \alpha & \gamma \cr \beta & \lambda },$
with $ \beta = 2 \gamma.$
Using the fact that $ \det (G)= \pm 1,$ we get
$ \beta ^2 \equiv \pm 2 \Mod\alpha.$
\EPf
\section{The harmonic knots $\mathbf{\H(a,b,c)}$}\label{harmonic}
In this paragraph we shall study Chebyshev knots with $ \phi= 0.$
Comstock (1897) found the number of crossing points of the
harmonic  curve parametrized by $x=T_a(t), y=T_b(t), z=T_c(t).$
In particular, he proved that this curve is non-singular if and
only if $ a,b,c$ are pairwise coprime integers (\cite{Com}). Such
curves will be named harmonic knots  $\H(a, b,c)$.
\pn
We shall need the following result proved in \cite{JP,KP3}
\begin{proposition}
Let $a $ and $b$ be coprime integers.  The $\frac 12 (a-1)(b-1)$
double points of the Chebyshev curve
$x= T_a(t), y= T_b(t)$ are obtained for the parameter pairs
$$
t= \cos \bigl( \Frac ka + \Frac hb \bigr) \pi,  \
s = \cos \bigl( \Frac ka - \Frac hb \bigr) \pi,
$$
where $h,k$ are positive integers such that
$ \Frac ka + \Frac hb < 1.$
\end{proposition}
Using the symmetries of Chebyshev polynomials, we see
that this set of parameters is symmetrical about the origin.
We will write  $ x \sim y$ when  $\sign x = \sign y.$
We shall need the following result proved in \cite{KP3}.
\begin{lemma}\label{sign}
Let $\H(a,b,c)$ be the harmonic knot: $x=T_a(t),\, y=T_b(t),\, z=T_c(t)$.
A crossing point of parameter
$
t = \cos \left (\Frac ka + \Frac hb\right ) \pi,  \
$
is a right twist if and only if
$$D = \Bigl( z(t)-z(s) \Bigr) x'(t) y'(t) >0$$
where
$$
z(t)-z(s) = T_c(t)-T_c(s) =
-2 \sin \Bigl( \Frac {ch}{b} \pi \Bigr) \sin \Bigl( \Frac {ck}a \pi \Bigr).
$$
and
$$
x'(t) y'(t) \sim (-1)^{h+k} \sin \Bigl( \Frac {ah}b \pi \Bigr)
\sin \Bigl( \Frac {bk}a \pi \Bigr).
$$
\end{lemma}
From this lemma we immediately deduce
\begin{corollary}\label{cprime}
Let $a,b,c$ be coprime integers. Suppose that the integer $c'$ verifies
$ c' \equiv c  \Mod{2a} $ and $ c' \equiv -c \Mod{2b}.$
Then the knot $\H(a,b,c')$ is the mirror image of $\H( a,b,c).$
\end{corollary}
\Pf
At each crossing point  we have
$
T_{c'}(t) - T_{c'}(s) = - \Bigl( T_c(t) - T_c (s)  \Bigr).
$
\EPf
\begin{corollary}
Let $a,b,c$ be coprime integers. Suppose that the integer $c$ is of
the form $c= \lambda a + \mu b$ with $\lambda, \mu >0$.
Then there exists $c'< c $
such that
$\H( a,b,c)= \overline{\H} (a,b,c') $
\end{corollary}
\Pf
Let $c'=\abs{\lambda a - \mu b}.$
The result follows immediately from corollary \ref{cprime}
\EPf
\pn
In a recent paper, G. and J. Freudenburg  have proved the following stronger result.
{\em There is a polynomial automorphism $ \Phi$
of $ \RR^3$ such that $ \Phi (\H(a,b,c)) = \H(a,b,c').$}
They also conjectured that the knots
$\H(a,b,c), \  a<b<c,  \  c \neq \lambda a + \mu b, \  \lambda, \mu >0 $
are  different knots (\cite{FF}, Conjecture 6.2).
%%%%%%%%%%%%%%%%%%%%%%%%%%%%%%%%%%%%%%%%%%%%%%%%%%%%%%%%%%%%%%%%%%%%%%%%%%%%%%%%%
\subsection{The harmonic knots ${\mathbf{\H(3,b,c)}}$ }
The following result is the main step in the classification of the
harmonic knots $\H( 3,b,c)$.
\begin{theorem}\label{h3}
Let $b=3n+1, \ c= 2b-3 \lambda, \, (\lambda, b)=1.$
The Schubert fraction of the knot $\H (3, b,c )$
is
$$ \Frac {\alpha}{\beta} =
[e_{1}, e_{2}, \ldots, e_{3n}],
\hbox{ where } e_k = \sign {\sin k \theta} \hbox{ and } \theta = \Frac {\lambda}b  \pi. $$
If $0< \lambda < \Frac {b}{2}$, its crossing number is
$N= b -  \lambda = \Frac {b+c}3,$
and we  have $\beta^2 \equiv \pm 1 \Mod \alpha.$
\end{theorem}
\Pf Will be given in section \ref{proofs}, p. \pageref{h3proof}.
\begin{corollary}\label{hh3}
The knots $\H(3,b,c)$ where $\Frac c2 <b<2c, \  b  \equiv 1 \Mod 3, \  c \equiv 2 \Mod 3$
are different knots (even up to mirroring). Their crossing number is given by $b+c=3N.$
\end{corollary}
\Pf
Let $K=\H(3,b,c)$ and $\Frac \alpha \beta>1$
be its biregular Schubert fraction given by Theorem \ref{h3}.
From Prop \ref{ell}, $\min(b,c)$ is
the minimum length of any Chebyshev diagram of  $K$ and
$\max(b,c)=3N-\min(b,c)$. The pair $(b,c)$ is uniquely determined.
\EPf

The following result gives the classification of harmonic knots $\H(3,b,c)$.
\begin{theorem}\label{h3bc}

Let $K= \H(3,b,c).$
 There exists  a
 unique  pair $(b',c')$
such that (up to mirror symmetry)
$$
K = \H(3,b',c'),   \  b'<c'<2b', \  b'+c' \equiv 0 \Mod 3.
$$
The crossing number of $K$ is $\frac 13 (b'+c'),$
its fractions $ \Frac \alpha \beta $ are such that
$ \beta^2 \equiv \pm 1 \Mod \alpha.$
Furthermore, there is an algorithm to find the pair $(b',c').$
\end{theorem}
\Pf
Let $K= \H(a,b,c)$
We will show that if the pair $(b,c)$ does not satisfy the condition of the theorem,
then it is possible to reduce it.

If $c<b$ we consider $\H(3,c,b)=\overline \H(3,b,c)$.

If $b \equiv c \Mod 3,$ we have $ c= b + 3 \mu, \  \mu>0.$
Let $c'= \abs{b-3\mu}$.
We have $c'\equiv \pm c \Mod {2b}$ and  $c' \equiv\mp c \Mod 6.$
By Lemma \ref{cprime}, we see that $ K= \overline \H(3,b,c')$ and we get a smaller
 pair.

If $b \not\equiv c \Mod 3$ and $c>2b,$ we have $c=2b + 3\mu, \ \mu>0.$
Let $c'= \abs{2b-3 \mu}$.
Similarly, we get  $ K= \overline\H(3,b,c')$.
This completes the proof of  existence.
This uniqueness is a direct consequence of  Corollary \ref{hh3}.
\EPf
%%%%%%%%%%%%%%%%%%%%%%%%%%%%%%%%%%%%%%%%%%%%%%%%%%%%%%%%%%%
\begin{remark}
Our theorem gives a positive answer to the Freudenburg conjecture for $a=3.$
\end{remark}
%\begin{remark}
%We thus see that the number of harmonic knots $\H(3,b,c)$
%of a given crossing number $N$ is very small.
%It is bounded by $\frac 12 \phi(N) \leq
%\frac 12 (N-1).$
%\end{remark}
\subsection*{Examples}
As  applications of Proposition \ref{bireg},
 let us deduce the following results (already   in \cite{KP3})
\begin{corollary}
The harmonic knot $ \H (3, 3n+2, 3n+1) $ is the torus knot $ T(2, 2n+1).$
\end{corollary}
\Pf
The harmonic knot $ K= \H (3, 3n+1, 3n+2 )$ is obtained for
$ b= 3n+1,$ $c= 2b -3\lambda, \  \lambda= n,$
%From Proposition \ref{bireg} its crossing number is $2n+1$.
%We have
 $ \theta = \Frac n {3n+1} \pi.$
If $j=1,2,$ or $ 3,$ and $ k=0,\ldots, n-1$
we have
$(3k+j) \theta = k \pi + \Frac {jk-n}{3n+1},$
hence $ \sign{ \sin (3k+j) \theta } = (-1)^k $, so that the Schubert fraction of $K$ is
$$
[1,1,1, -1,-1,-1, \ldots, (-1)^{n+1}, (-1)^{n+1}, (-1)^{n+1} ] = \Frac{2n+1} {2n}
\thickapprox -(2n+1).
$$
We see that $K$ is the mirror image of $ T(2,2n+1),$ which completes the proof.
\EPf
\pn
It is possible to parameterize the knot $T(2, 2n+1)$ by polynomials of the same degrees
and a diagram with only $2n+1$ crossings (\cite{KP2}).
However, our Chebyshev parametrizations are easier to visualize.
We conjecture that these degrees are minimal
(see also \cite{RS,FF,KP1}).

\begin{corollary}
The harmonic knot $\H(3,b,2b-3)$ ($b\not \equiv 0 \Mod 3$)  is alternate and has
crossing number $b-1$.
\end{corollary}
\Pf
For this knot we have $ \lambda=1, \  \theta = \Frac \pi b.$
The Schubert fraction is given by the continued fraction of length $b-1$:
[1,1, \ldots,1] = $\Frac{F_{b}}{F_{b-1}}$ where $F_n$ are the Fibonacci numbers
($F_0=0, F_1=1, \ldots$).
J. C. Turner named these  knots  Fibonacci knots (\cite{Tu}).
\EPf
\pn
The two previous examples describe infinite families of harmonic knots. They have
a Schubert fraction $\Frac\alpha\beta$ with $\beta^2 \equiv  1 \Mod\alpha$ (torus knots) or
with $\beta^2 \equiv - 1\Mod\alpha$ (Fibonacci knots with odd $b$). There
is also an infinite number of
two-bridge knots with $\beta^2= \pm 1 \Mod\alpha$ that are not harmonic.
\begin{proposition}
The knots (or links)
$K_n = S(\Frac{5F_{n+1}}{F_{n+1}+F_{n-1}}),\, n>1$ are not harmonic knots
$\H(3,b,c)$. Their crossing number is $n+4$ and we have
$(F_{n+1}+F_{n-1})^2 \equiv (-1)^{n+1} \Mod{5F_{n+1}} $.
\end{proposition}
\Pf
Using the fact that $P^n = \cmatrix{F_{n+1}&F_n\cr F_n&F_{n-1}}$ we deduce
that
$$
P M P^n M P = \cmatrix{5F_{n+1}&F_{n+1}+F_{n-1}\cr F_{n+1}+F_{n-1}&F_{n-1}}.
$$
Taking  determinants, we get $(F_{n+1}+F_{n-1})^2 \equiv (-1)^{n+1} \Mod{5F_{n+1}} $.
We also have
$$
\Frac{5F_{n+1}}{F_{n+1}+F_{n-1}} =
P M P^n M P(\infty)=
[1,1,\underbrace{-1,\ldots,-1}_{n+2},1,1].
$$
Since $ n+2 \ge 4,$ it cannot be of the form
  $[\sign{\sin \theta},\sign{\sin 2\theta}, \ldots, \sign{\sin k\theta}]$.

If $n\equiv 2\Mod 3$, $K_n$ is a two-component link.

If $n\equiv 1\Mod 6$ or $n\equiv 3\Mod 6$, the Schubert fraction $\Frac\alpha\beta$
satisfies $\beta^2\equiv 1\Mod\alpha$.

If $n\equiv 0\Mod 6$ or $n\equiv 4\Mod 6$, $K_n$ is amphicheiral.
\EPf
\subsection{The harmonic knots ${\mathbf{\H(4,b,c)}}.$}
The following result will allow us to  classify the harmonic knots of the form
 $\H(4,b,c).$
\begin{theorem}\label{h4}
Let $b, c$ be odd integers such that $ b \not \equiv c  \Mod 4.$
The Schubert fraction of the knot $K= \H(4,b,c)$ is given by the
 continued fraction
$$
\Frac \alpha\beta =
[ e_1, 2e_2, e_3, 2e_4, \ldots, e_{b-2}, 2e_{b-1} ], \,
e_j= -\sign{\sin (b-j) \theta}, \,
\theta = \Frac {3b-c} {4b}  \pi.$$
If $b< c < 3b,$ this fraction is biregular, the crossing number of $K$ is
$ N= \Frac { 3b+c-2}4,$ and $\beta ^2  \equiv \pm 2 \Mod 3.$
\end{theorem}
\Pf Will be given in section \ref{proofs}, p. \pageref{h4proof}.\EPf
\pn
We are now able to classify  the harmonic knots of the form $\H(4,b,c)$.
\begin{theorem}\label{h4bc}\\
Let $K= \H(4,b,c).$
There is a unique pair
$(b',c')$
such that (up to mirroring)
$$K= \H(4,b',c'), \    b'<c'<3b', \  b' \not\equiv c' \Mod 4.$$
The crossing number of $K$ is $\frac 14 (3b'+c'-2). $
$K$ has a   Schubert fraction
$\Frac \alpha \beta $  such that\\
$ \beta^2 \equiv \pm 2 \Mod \alpha.$
Furthermore, there is an algorithm to find $(b', c').$
\end{theorem}
\Pf
Let us first prove the uniqueness of this pair.
Let $K= \H(4,b,c)$ with $ b<c<3b, \ c \not\equiv b \Mod 4.$
By  theorem \ref{h4}, $K$ has a Schubert fraction
$ \Frac \alpha \beta = \pm [ 1, \pm2, \ldots, \pm 1, \pm 2] $
of length $b-1$ with $\beta$ even, $ -\alpha < \beta < \alpha,$ and
$ \beta ^2 \equiv \pm 2 \Mod \alpha.$
The other fraction of $K$ is
$\Frac \alpha{\beta'}$, where $\beta'$ is even and
 $ \abs{\beta'}<\alpha $ and $ \beta   \beta ' \equiv 1 \Mod \alpha.$

If $ \alpha > 3,$ we cannot have $ {\beta'}^2 \equiv \pm 2 \Mod \alpha.$
We deduce that $b$ is uniquely determined by $K$:
$b=\ell ( \Frac \alpha \beta ) +1 $ where $\Frac \alpha \beta$ is
a Schubert fraction of $K$ such that $\beta^2 \equiv \pm 2 \Mod \alpha.$
Since we also have $3 b+c -2= 4 \cn(K),$
we conclude that $(b,c)$ is uniquely determined by $K.$
\pn
Now, let us prove the existence of the pair $(b',c').$
Let $K= \H(4,b,c),  \   b<c$.
We  have only to  show that if the pair $(b,c)$
does not satisfy
the condition of the theorem, then it is possible to reduce it.
\pn
If $ c \equiv b \Mod 4, $ then $c= b + 4 \mu, \  \mu >0.$
Let $c'= \abs{b - 4 \mu}$.
Then  $K= \overline \H( 4, b, c' )$, and the pair $(b,c')$ is smaller than $(b,c).$
\pn
If $ c \not\equiv b \Mod 4 $ and $c>3b,$
we have $c= 3b + 4 \mu, \  \mu >0.$
Let  $c' = \abs{3b-4 \mu}$.
We have, $K= \overline\H(4,b,c')$ with $(b,c')$ smaller than $(b,c).$
This completes the  proof.
\EPf
\begin{remark}
Our theorem gives a positive answer to the Freudenburg conjecture for $a=4.$
\end{remark}
We also see that the only knot belonging to the two families
of  knots $ \H(3,b,c)$ and $\H(4,b,c)$
is the trefoil $ \H(3, 4,5 )= \overline \H(4,3,5).$
\begin{example}[$\mathbf{\H(4,2k-1,2k+1)}$]\\
From Theorem \ref{h4}, we know that $\H(4,2k-1,2k+1)$ has crossing number
$2k-1$. Using this theorem, the Conway sequence of this knot is
$
[e_1,2e_2, \ldots, e_{2k-3},2e_{2k-2}]$, where
$$e_j = -\sign{\sin((2k-1-j)\Frac{(k-1)\pi}{2k-1})}
= (-1)^{k+1} \sign{\sin(\Frac{j(k-1)\pi}{2k-1})} = (-1)^{k+\pent{j+1}2}.
$$
We deduce that the Schubert fraction of $\H(4,2k-1,2k+1)$ is
$$
\Frac{\alpha_k}{\beta_k} =
(-1)^{k+1} [1,2,-1,-2,1,2,\ldots, (-1)^{k},2(-1)^{k}] = (-1)^{k+1} (AS)^{k-2}A (\infty).
$$
Using recurrence formula in Lemma \ref{nn}, we deduce that
$$
\begin{array}{rcl}
\alpha_2 = 3, & \alpha_3=7, & \alpha_{k+2} = 2 \alpha_{k+1}+\alpha_k\\
  \beta_2 = -2, & \beta_3=4, & \abs {\beta_{k+2}} = 2 \abs {\beta_{k+1}}+\abs {\beta_k}.
\end{array}
$$
Let us consider the homography
$ G(x)= [2,x],$
and its matrix
$ G = \cmatrix { 2 & 1 \cr 1 & 0 }$.
Let the sequence $u_k$ be defined by
$$  G^k = \cmatrix { u_{k+1} & u_k \cr u_k & u_{k-1}}.$$
The sequence $u_k$ verifies the same recurrence formula
$ u_{k+2}= 2 u_{k+1} + u_k. $ We deduce
$\alpha _k = u_k + u_{k-1}, \abs { \beta _k } = 2 u_{k-1}.$
We also have
$$
P^2 G^{k-2} P ( \infty ) =
\cmatrix { 2 & 1 \cr 1 & 1 } \,
\cmatrix { u_{k-1} & u _ {k-2} \cr u_{k-2} & u _{k-3} } \,
\cmatrix {1 \cr 1 }
= \cmatrix { u_k + u_{k-1} \cr 2 u_k }
= \Frac {\alpha _k }{ \abs { \beta _k } }.
$$
Finally, we get  $r_k = (-1)^{k+1} [1,1,\underbrace{2, \ldots 2}_{k-2},1]$.
\end{example}
\begin{example}[The twist knots]\label{twh4}
The  knots $\cT_n$  are not harmonic knots $\H(4,b,c) $ for $n>3.$
\end{example}
\Pf
The Schubert fractions of $\cT_n= S( n+ \Frac 12 ) $ with an even denominator
are
$ \Frac {2n+1}2,$ and  $ \Frac{2n+1}{-n} $ or $ \Frac {2n+1} {n+1} $ according
to the parity of $n.$
The only such fractions verifying $ \beta ^2 \equiv \pm 2 \Mod \alpha $ are
$ \Frac 32, \, \Frac 74, \, \Frac 94.$
The first two are the Schubert fractions of the trefoil and the $\overline{5}_2$ knot,
which are harmonic for $a=4.$
We have only to study the case of $ 6_1 = S ( \Frac 9{4} ).$
We have
$ \Frac 9{4} =  [ 1,2, -1,2, 1,-2, 1,2 ].$
Since this fraction is not biregular, we see that $6_1$ is not of the form
$\H(4,b,c).$
\EPf
\pn
But there also exist  infinitely many rational knots whose Schubert fractions
$\Frac\alpha\beta$ satisfy $\beta^2\equiv -2 \Mod\alpha$  that are not
harmonic for $a=4$.
\begin{proposition}
The knots $S(n+\Frac{1}{2n})$ are not harmonic knots $\H(4,b,c)$ for $n>1$.
Their crossing number is $3n$ and their Schubert fractions $\Frac\alpha\beta = \Frac{2n^2+1}{2n}$
satisfy $\beta^2\equiv -2 \Mod \alpha$.
\end{proposition}
\Pf
If $n=2k$, we deduce
from $(ASB)^k(x) = 2k+x$ and $(BSA)^k (\infty)= \Frac{1}{4k}$, that
$$n+\Frac{1}{2n} = (ASB)^k (BSA)^k(\infty).$$
If $n=2k+1$, we use (see the torus knots, example \ref{tnh4})
$\Frac{2n+1}{2n} = A(BSA)^k(\infty)$, so
$$n+\Frac{1}{2n} = n-1 + \Frac{2n+1}{2n} = (ASB)^k A(BSA)^k(\infty).$$
For $n>1$ these continued fractions  are not biregular, and since
$ \beta ^2 \equiv -2 \Mod \alpha$, they do
not correspond to harmonic knots $\H(4,b,c)$.
\EPf
\section{Chebyshev diagrams of rational  knots}\label{diagrams}
\begin{definition} We say that a knot in $\RR ^3  \subset  {\SS}^3$ has
a Chebyshev diagram $\cC(a,b)$, if $a$ and $b$ are coprime and the Chebyshev curve
$$\cC(a,b): x=T_a(t); \  y=T_b(t)$$ is the projection of some knot which is
isotopic to $K$.
\end{definition}
\subsection{Chebyshev diagrams with $a=3$}
Using the previous results of our paper (Proposition \ref{ell}) we have
\begin{theorem}
Let $K$ be a two-bridge knot with crossing number $N$.
There is an algorithm to determine the smallest $b$ such that
$K$ has a Chebyshev diagram $\cC(3,b)$ with $N < b < \Frac 32 N$.
\end{theorem}
\Pf
Let $\Frac \alpha \beta > 1$ and $ \Frac \alpha{\alpha - \beta}>1 $ be
Schubert fractions of $K$ and $\overline{K}.$ By Proposition \ref{ell},
$b=\min\Bigl( \ell (\Frac \alpha \beta),\ell (\Frac \alpha {\alpha - \beta})\Bigr)+1$
has the required property.
\EPf
\begin{definition}
Let ${\cal D}(K)$ be a diagram of a knot having crossing points corresponding to
the parameters  $ t_1, \ldots, t_{2m} $.
The Gauss sequence of ${\cal D}(K)$ is defined by
$ g_k = 1 $ if $t_k$ corresponds to an overpass, and
$ g_k = -1$ if $t_k$ corresponds to an underpass.
\end{definition}
\begin{theorem}\label{gauss3}
Let $ K$ be a two-bridge knot of crossing number $N.$
Let $x= T_3(t),\ y =T_b(t)$ be the minimal Chebyshev diagram of $K$.
Let $c$ denote the number of sign changes in the corresponding Gauss sequence.
Then we have
$$ b+c=3N.$$
\end{theorem}
\Pf Let $s$ be the number of sign changes in the Conway normal form of $K.$
By Proposition \ref{bireg} we have $ N= b-1-s.$
From this we deduce that our condition is equivalent to $3s+c = 2b-3.$
Let us prove this assertion by induction on $s.$
If $s=0$ then the diagram of $K$ is alternate,
and we deduce $c= 2(b-1) -1=2b-3.$

Let $C(e_1,e_2, \ldots, e_{b-1}) $ be the Conway normal form of $K.$ We may suppose
$e_1=1$. We shall denote by $M_1, \ldots, M_{b-1}$ the crossing points of the diagram, and
by $ x_1< x_2 < \cdots < x_{b-1} $ their abscissae.
Let $e_k$ be the first negative coefficient in this form.
By the regularity of the sequence we get $e_{k+1}<0,$ and $ 3 \le k \le b-1.$

Let us consider the knot $K'$ defined by its Conway normal form
$$ K' = C(e_1,e_2, \ldots, e_{k-1}, -e_k, -e_{k+1}, \ldots, - e_{b-1}).$$
We see that  the number of sign changes in the Conway sequence of $K'$ is
$s'= s-1.$ By induction, we get for the knot $K'$:  $3s'+c'= 2b-3.$

The  plane curve $x=T_3(t), \  y= T_b(t) $ is the union of three
arcs where $x(t)$ is monotonic.
Let $ \Gamma$ be one of these arcs.
$\Gamma$ contains (at least) one point $M_k$ or $M_{k+1}.$
\psfrag{g1}{\large$\mathbf{\Gamma_1}$}%
\psfrag{g2}{\large$\mathbf{\Gamma_2}$}%
\psfrag{g3}{\large$\mathbf{\Gamma_3}$}%
\psfrag{k1}{$M_{k-1}$}%
\psfrag{k2}{$M_k$}%
\psfrag{k3}{$M_{k+1}$}%
\begin{figure}[th]
\begin{center}
\begin{tabular}{ccc}
{\scalebox{.7}{\includegraphics{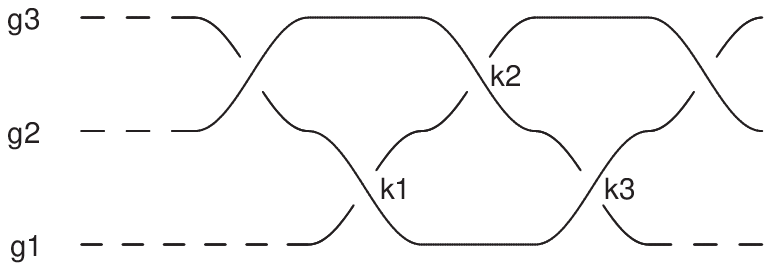}}}&&
{\scalebox{.7}{\includegraphics{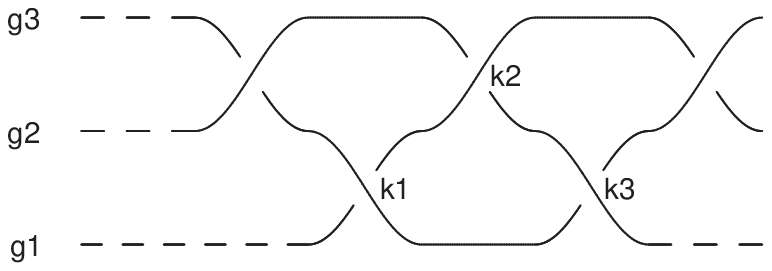}}}\\
$K$&\hspace{2cm}&$K'$
\end{tabular}
\end{center}
\caption{The modification of Gauss sequences}
\label{gaussf}
\end{figure}
Let $j$ be the first integer in $\{ k, k+1 \}$ such that $M_j$ is on $\Gamma$,
and let $j_- < j$ be the greatest integer such that $ M_{j_-} \in \Gamma.$
In figure \ref{gaussf}, we have
for $\Gamma_1$: $j=k, j^- = k-1$,
for $\Gamma_2$: $j=k, j^- = k-2$,
for $\Gamma_3$: $j=k+1, j^- = k-1.$

On each arc $\Gamma$, there is a sign change in the Gauss sequence iff the corresponding Conway
signs are equal.
Then, since the Conway signs $ s(M_{j_-} )$ and $s(M_j)$ are different, we see that
the corresponding Gauss signs are equal.
Now, consider the modifications in the Gauss sequences when we transform
$K$ into $K'$.
Since the the Conway signs $s( M_h), \, h\ge k $ are changed,
we see that we get one more sign change on  every arc $\Gamma$.
Thus the number of sign changes in the Gauss sequence of $K'$ is
$c' = c+3.$
We get $3s+c= 3(s'+1) + c'-3= 3s'+c' = 2b-3,$ which completes our induction proof.
\EPf
\begin{corollary}\label{degc3}
Let $K$ be a two-bridge knot with crossing number $N$. Then there exist $b,c$, $b+c=3N$,
and an  polynomial $C$ of degree $c$ such that the knot
$x=T_3(t), \, y=T_b(t), \, z=C(t)$ is isotopic to $K$.

If $K$ is amphicheiral, then  $b$ is odd, and the polynomial $C(t)$
can be chosen odd.
\end{corollary}
\Pf Let $b=n+1$ be the smallest integer  such that $K$ has
a Chebyshev diagram $x=T_3(t), \, y= T_b(t)$.
By our theorem \ref{gauss3}, the Gauss sequence $(g(t_1), \ldots, g(t_{2n}))$
of this diagram has $c=3N-b$ sign changes. We choose $C$ such that
$C(t_i)g(t_i)>0$ and we can realize it by choosing the roots of
$C$ to be $\frac 12(t_{i}+t_{i+1})$ where $g(t_i)g(t_{i+1})<0$.

If $K$ is amphicheiral, then  $b$ is odd and the Conway form is palindromic
by Proposition \ref{palin}.
Then our Chebyshev diagram is symmetrical about the origin. We see that the Gauss
sequence is odd: $g(t_h) = -g(-t_h)$.
This implies that the polynomial $C(t)$ is odd.
\EPf
\begin{remark}
This corollary gives a simple proof of a famous theorem
of Hartley and Kawauchi:
{\em every amphicheiral rational knot is strongly negative amphicheiral}
(\cite{HK,Kaw}).
\end{remark}
\begin{example}[The knot $6_1$]\label{6_1-3}\\
The knot $\overline{6}_1= S(\Frac 92)$ is not harmonic with $a=4.$  It is not even harmonic with $a=3$
because $2^2 \not \equiv \pm 1 \Mod 9$.
Its crossing number is $6$.
In the example \ref{9297}, we get $\ell(\Frac 92)=9,\, \ell(\Frac 97)=7$. $b=8$ is the minimal value for which
$x=T_3(t), \, y=T_8(t)$ is a Chebyshev diagram for $\overline{6}_1$.
The Gauss sequence associated to the Conway form
$\overline{6}_1 = C(-1,-1,-1,1,1,1,1)$ has exactly 10 sign changes. It is precisely
$$[1, -1, -1, 1, -1, 1, -1, -1, 1, -1, 1, 1, -1, 1].$$
We can check that
$$
x=T_3(t),\, y=T_8(t),\, z=
 \left( 8\,t+7 \right)  \left( 5\,t-4 \right)  \left( 15\,{t}^{2}-14
 \right)  \left( 2\,{t}^{2}-1 \right)  \left( 3\,{t}^{2}-1 \right)
 \left( 15\,{t}^{2}-1 \right)
%\left( 9\,{t}^{2}-4 \right)  \left( 900\,{t}^{2}-841 \right)
%\left( 4\,{t}^{2}-1 \right)  \left( 16\,{t}^{2}-1 \right)  \left(40
%\,{t}^{2}-3\,t-28 \right)
$$
is a parametrization of $\overline{6}_1$ of degree $(3,8,10)$.  In \cite{KP3} we gave
the Chebyshev parametrization $6_1= \cC(3,8,10,\frac{1}{100})$. We will give another
parametrization in example \ref{6_1-4}.
\begin{figure}[th]
\begin{center}
{\scalebox{.7}{\includegraphics{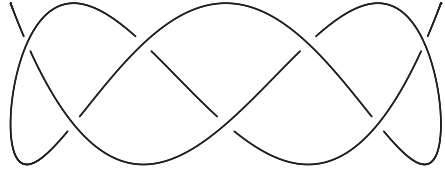}}}
\end{center}
\caption{The knot $6_1$}
\label{61_3}
\end{figure}
\end{example}
\subsection{Chebyshev diagrams with $a=4$}
It is also possible to get Chebyshev diagrams of the form $\cC(4,b)$.
The following result is analogous to the Theorem \ref{gauss3}.
\begin{theorem}\label{gauss4}
Let $K$ be a two-bridge knot of crossing number $N$ and Schubert fraction
$\Frac \alpha \beta$, $\beta$ even.
Let $\Frac \alpha \beta = \pm [1,\pm2, \ldots,\pm1, \pm 2 ]$ be a continued
fraction expansion of minimal length $b-1,$ and
$\sigma$ be the number of islets (subsequences of the form $ \pm ( 2, -1, 2 )$) in
this expansion.
Let $x= T_4(t),\, y = T_b(t) $ be the corresponding Chebyshev diagram of $K$
and $c$ be the number of sign changes in the corresponding Gauss sequence of $K$.
Then we have
$$
3b+c-2 = 4N + 12 \sigma.
$$
\end{theorem}
\Pf Let $s$ be
the number of sign changes in the given continued fraction. Since
$N= \frac 32 (b-1) -s -2 \sigma $ the formula is equivalent to
$3b+c-2 = 6( b-1) -4s -8 \sigma + 12 \sigma, $ that is
$$
4 s +c - 4 \sigma  = 3b -4.
$$
We shall prove this formula by induction on $s.$

If $s=0,$ the knot is alternate. We have  $ c = 3b-4, \ \sigma
=0,$ and the formula is true.

We shall need precise notations. Let $ M_1, M_2,N_2, M_3, M_4,
N_4,\ldots, M_{b-1}, N_{b-1} $ be the crossing points where $
x(M_k)=x_k, \  x(N_{2k})= x_{2k}, $ and $ x_1< x_2< \cdots <
x_{b-1}.$

The plane curve $x= T_4(t), \ y= T_b(t) $ is the union of four
arcs where $x(t)$ is monotonic (see the following figures).
On each arc there is
a sign change in the Gauss sequences iff the corresponding Conway
signs are equal. Let $ C( e_1,2e_2, \ldots, e_{b-2}, 2 e_{b-1} ),
\  e_i= \pm1$ be the Conway form of $K.$ Let $k$ be the first
integer such that $ e_{k-1} e_k <0.$ We have three cases to
consider.
\pn
{\bf $\mathbf{k}$ is odd, and $\mathbf{e_k e_{k+1} <0}$.}\\
In this case, $ ( 2e_{k-1}, e_k, 2e_{k+1})= \pm (2, -1, 2) $ is an
islet. Let us consider the knot $K'$ obtained by changing only the
sign of $e_k.$ The number of sign changes in the Conway sequence
of $K'$ is $s'=s-2.$ By induction we get for $K':$ $4s' +c' -4
\sigma' = 3b-4.$ The number of islets of $K'$ is $ \sigma ' =
\sigma -1.$ Let us look the modification of Gauss sequences.
There are only two arcs containing $ M_k.$ On each of these arcs
there is no sign change in the Gauss sequence of $K,$ and then
there are two sign changes in the Gauss sequence of $K'.$
Consequently, $c'= c+4.$
\psfrag{g1}{\large$\mathbf{\Gamma_1}$}%
\psfrag{g2}{\large$\mathbf{\Gamma_2}$}%
\psfrag{g3}{\large$\mathbf{\Gamma_3}$}%
\psfrag{g4}{\large$\mathbf{\Gamma_4}$}%
\psfrag{a1}{$M_{k-1}$}%
\psfrag{c1}{$N_{k-1}$}%
\psfrag{b2}{$M_{k}$}%
\psfrag{a2}{$M_{k+1}$}%
\psfrag{c2}{$N_{k+1}$}%
\begin{figure}[th]
\begin{center}
\begin{tabular}{ccc}
{\scalebox{.7}{\includegraphics{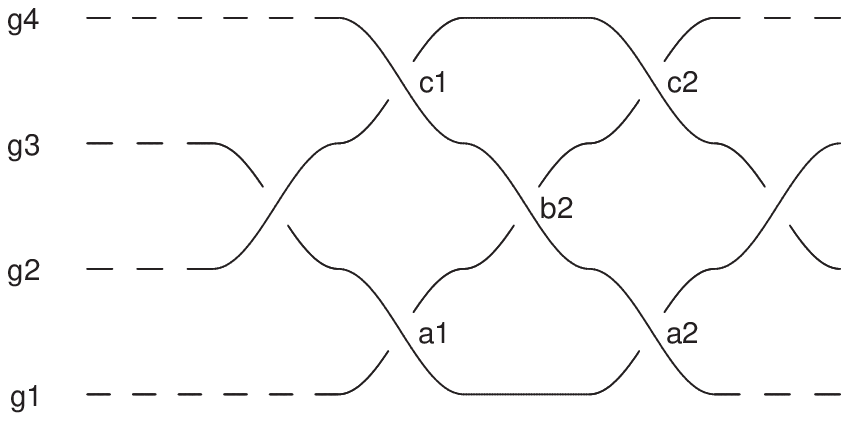}}}&&
{\scalebox{.7}{\includegraphics{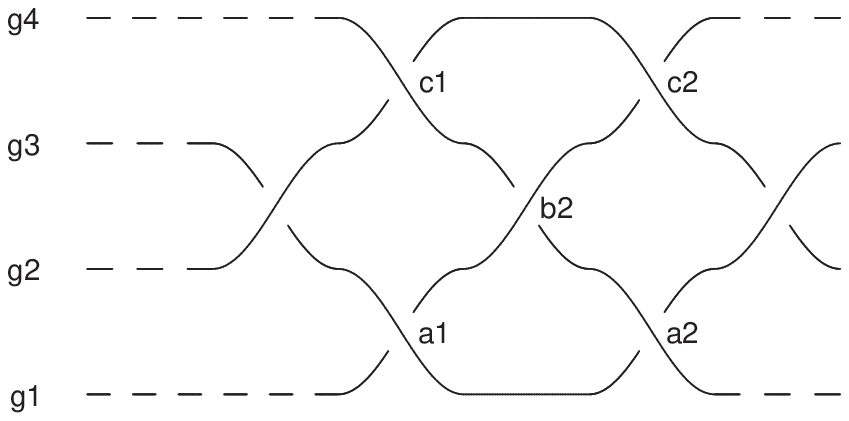}}}\\
$K$&\hspace{2cm}&$K'$
\end{tabular}
\end{center}
%\caption{The modification of Gauss sequences}
\label{dgauss41}
\end{figure}
Finally, we get $4s +c -4 \sigma =  4 ( s'+2) + (c'-4) -4( \sigma
' +1) = 4s' +c' -4 \sigma ' = 3b-4,$ which completes the proof in
this case.
\pn
{\bf $\mathbf{k}$ is even, and $\mathbf{e_k e_{k+1} <0}$.}\\
In this case there are two crossing points $ M_k$ and $N_k$ with $
x(M_k)=x(N_k)= x_k.$ Each arc contains one of these points. Let us
consider the knot $K'$ obtained by changing only the sign of
$e_k.$ The number of sign changes in the Conway sequence of $K'$
is $s' = s-2.$ By induction the formula is true for $K'.$
\psfrag{b0}{$M_{k-1}$}%
\psfrag{a1}{$M_{k}$}%
\psfrag{c1}{$N_{k}$}%
\psfrag{b2}{$M_{k+1}$}%
\psfrag{a2}{}%
\psfrag{c2}{}%
\begin{figure}[th]
\begin{center}
\begin{tabular}{ccc}
{\scalebox{.7}{\includegraphics{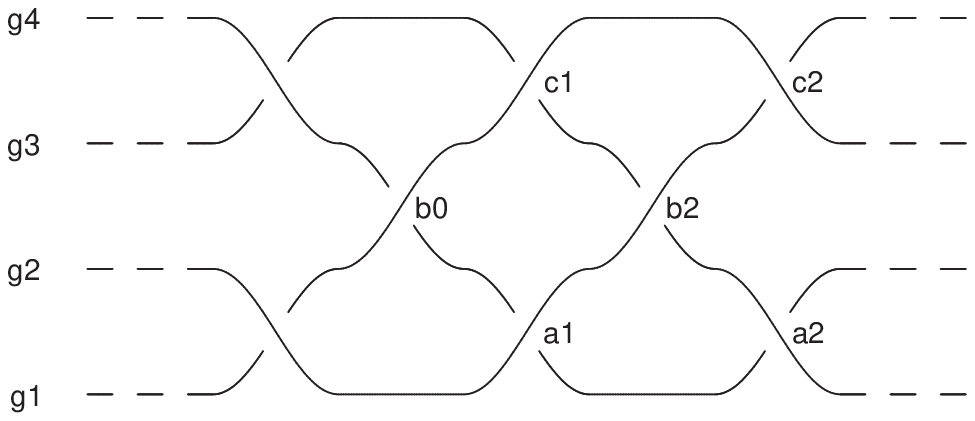}}}&&
{\scalebox{.7}{\includegraphics{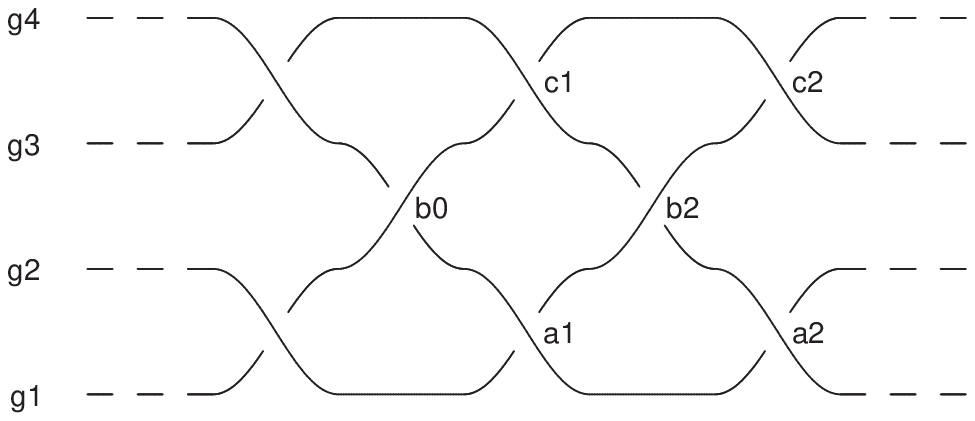}}}\\
$K$&\hspace{2cm}&$K'$
\end{tabular}
\end{center}
%\caption{The modification of Gauss sequences}
\label{dgauss42}
\end{figure}
On each arc the Gauss sequence gains two sign changes, so that
$c'=c+8.$ Since we have $\sigma ' = \sigma,$ we get
$
4 s + c - 4 \sigma = 4 ( s'+2) + (c'-8) -4 \sigma ' = 3b-4.$
\pn
{\bf The case $\mathbf{e_k e_{k+1} >0}$.}\\
\begin{figure}[th]
\begin{center}
\begin{tabular}{ccc}
{\scalebox{.7}{\includegraphics{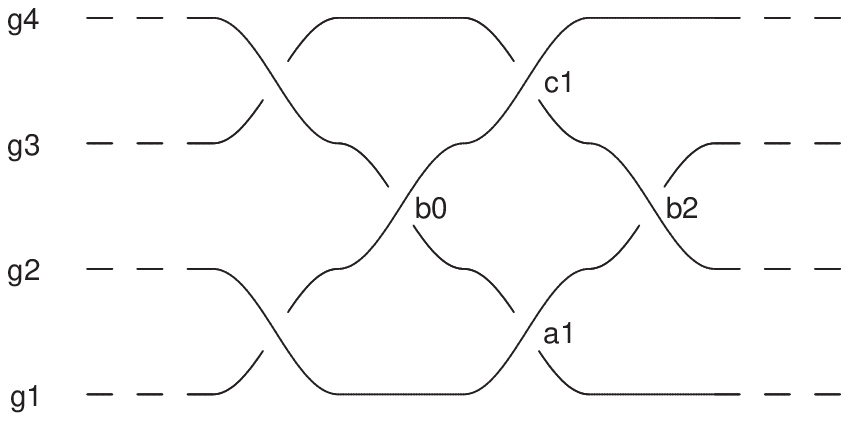}}}&&
{\scalebox{.7}{\includegraphics{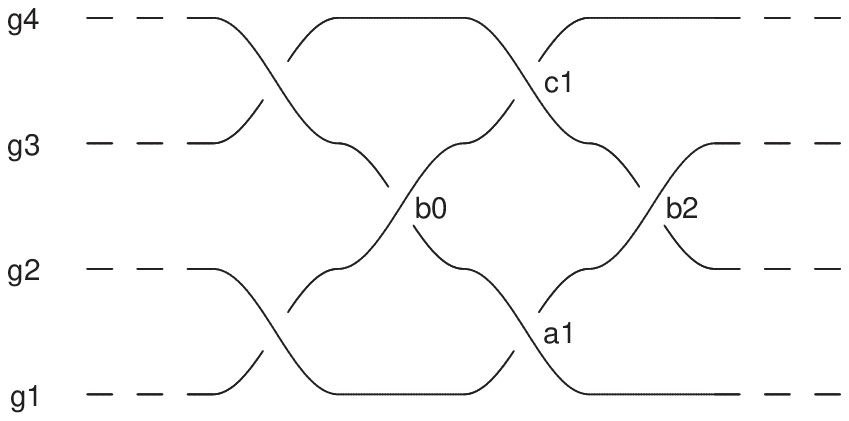}}}\\
$K$&\hspace{2cm}&$K'$
\end{tabular}
\end{center}
%\caption{The modification of Gauss sequences}
\label{dgauss43}
\end{figure}%
In this case we consider the knot $K'$ obtained by changing the
signs of $e_j, \  j \ge k.$ For $K'$ we have $ s' = s-1,$ and by
induction $ 4 s' + c' - 4 \sigma' = 3b-4.$\\
On each of the four arcs the Gauss sequence gains one sign change,
and then $c'=c+4.$ Since $ \sigma ' = \sigma,$ we conclude
$$4 s+c- 4 \sigma= 4 (s'+1) + (c'-4) - 4 \sigma ' =
4 s' + c' - 4 \sigma ' = 3b-4.$$ This completes the proof of the
last case.
\EPf
\begin{corollary}\label{degc4}
Let $K$ be a two-bridge knot of crossing number $N$ and Schubert fraction
$\Frac \alpha \beta$, $\beta$ even.
Let $\Frac \alpha \beta = \pm [1,\pm2, \ldots,\pm1, \pm 2 ]$ be a continued
fraction expansion of minimal length $b-1$ and
$\sigma$ be the number of islets (subsequences of the form $ \pm ( 2, -1, 2 )$) in
this expansion. There exists an odd polynomial $C(t)$ of degree $c$ such that
$3b+c-2 = 4N + 12 \sigma$ and such that the knot defined by
$x=T_4(t),\, y=T_b(t), \, z=C(t)$ is isotopic to $K$.
\end{corollary}
\Pf
This proof is similar to the proof of Corollary \ref{degc3}. We must bear
in mind that in this case, the Gauss sequence is odd:
$g(t_i) = -g(-t_i)$ and $t_{3(b-1)+1-i} = -t_i$.
\EPf
\pn
Corollary \ref{degc4} gives an algorithm to represent any rational knot as a
polynomial knot which is rotationally symmetric around  the $y$-axis.
This gives a strong evidence for the following classical result.
\begin{corollary}
Every rational knot is strongly invertible.
\end{corollary}
\begin{example}[The stevedore knot $\overline{6}_1$]\label{6_1-4}
The stevedore knot $6_1$ is $S(\Frac 92)=S(-\Frac 94)$.
We get $\Frac 94 = (ASB)(BSA)(\infty) = [1,2,-1,2,1,-2,1,2]$. We deduce that it can be parametrized
by $x=T_4(t),\ y= T_9(t), \ z=C(t)$ where $\deg C = 11$. We find
$$
C(t) = t \left( 34\,{t}^{2}-33 \right)  \left( 2\,{t}^{2}-1 \right)  \left( 3
\,{t}^{2}-1 \right)  \left( 4\,{t}^{2}-1 \right)  \left( 6\,{t}^{2}-1
\right).
$$
On the other hand, we find $S(\Frac 9{14}) = [1,-2,-1,-2,-1,2] = BASB(\infty)$. We find
that $6_1$ can also be represented by polynomials of degrees $(4,7,9)$:
$$
x=T_4(t),\ y= T_7(t), \ z=t \left( 10\,{t}^{2}-9 \right)  \left( 4\,{t}^{2}-3 \right)  \left( 4
\,{t}^{2}-1 \right)  \left( 6\,{t}^{2}-1 \right).
$$
\begin{figure}[th]
\begin{center}
\begin{tabular}{ccc}
{\scalebox{.7}{\includegraphics{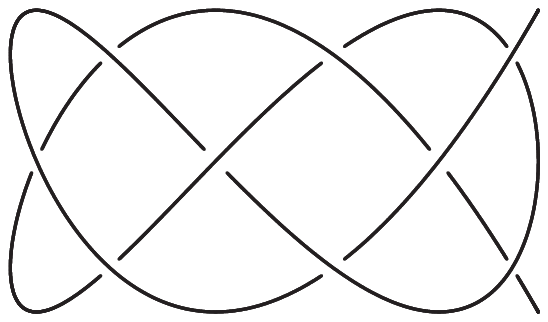}}}&&\\
%{\scalebox{.7}{\includegraphics{gauss43b.eps}}}\\
$\overline{6}_1 = S(\Frac{9}{14})$&&
\end{tabular}
\end{center}
%\caption{The modification of Gauss sequences}
\label{61_4}
\end{figure}%
\end{example}
\subsection{Examples}
In this section, we give several examples of polynomial parametrizations of rational knots
with Chebyshev diagrams
$\cC(3,b)$ and $\cC(4,b').$
%Corollary \ref{degc4} allows us to give explicit
%polynomial parametrization with Chebyshev diagrams $\cC(4,b)$.
\subsubsection*{Parametrizations of the torus knots}
The torus knot $T(2,N), \ N=2n+1$ is the harmonic knot
 $\overline\H(3,3n+1,3n+2)$.

The torus knot  $T(2, N)$ can be parametrized by
$x=T_4(t), \,  y= T_N (t), \, z= C(t),$ where $C(t)$ is an odd polynomial
of degree  $\deg (C) = N+2=2n+3.$
\pn
\Pf
We have already seen (example \ref{tnh4}) that $T(2,N)$ has a Chebyshev diagram $ x= T_4(t), y= T_N(t)$.
Since there is no islet in the corresponding continued fractions, we see
that the number of sign changes in the Gauss sequence is $c= N+2.$
By symmetry of the diagram, we can find an odd polynomial of degree
$c$ giving this diagram.
\EPf
\pn
In both cases, the  diagrams
 have the same number of crossing points : $\frac{3}2{(N-1)}=3n$.
\pn
As an example, we obtain for $\overline{T}(2,5)$:
$$
x = T_4(t),\, y=T_5(t), \, z =
t \left( 2\,{t}^{2}-1 \right)  \left( 3\,
{t}^{2}-1 \right) \left( 5\,{t}^{2}-4 \right).
$$
In this case the Chebyshev diagram has exactly 6 crossing points as it is
for $\H(3,7,8)$.
Note that we also obtain  $ T(2,5)=S(\Frac 56)$:
$$
x = T_4(t),\, y=T_7(t), \, z =t \left( 21\,{t}^{2}-20 \right)  \left( 4\,{t}^{2}-1 \right).
$$
We therefore obtain parametrizations of degrees $(4,5,7)$ or $(4,7,5)$.
\begin{figure}[th]
\begin{center}
\begin{tabular}{ccc}
{\scalebox{.7}{\includegraphics{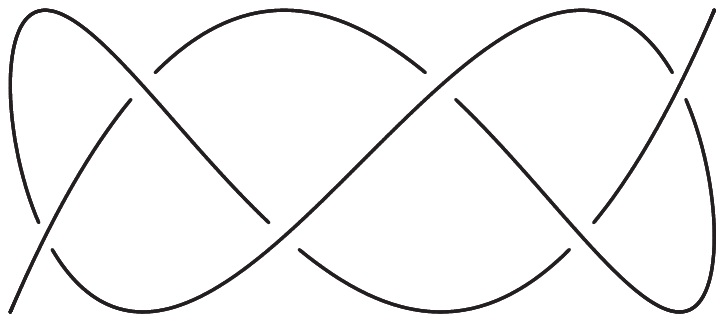}}}&
{\scalebox{.7}{\includegraphics{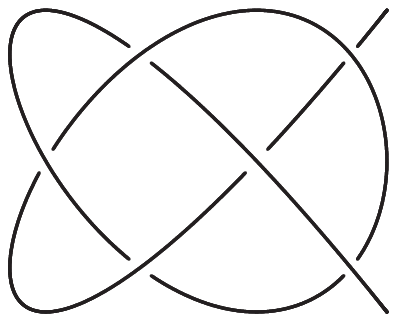}}}&
{\scalebox{.7}{\includegraphics{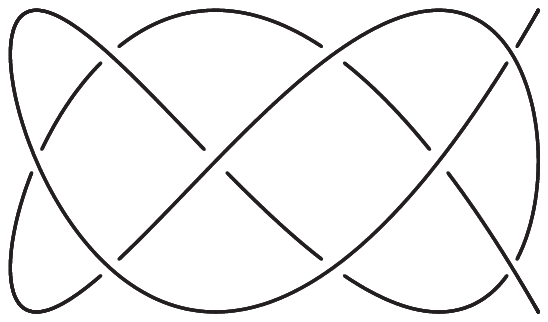}}}\\
$\H(3,7,8)$&$S(\Frac 54)$&$S(\Frac 56)$
\end{tabular}
\end{center}
\caption{Diagrams of the torus knot $\overline{T}(2,5)$ and its mirror image}
\label{K51}
\end{figure}
\subsubsection*{Parametrizations of the twist knots}
The twist knot $\cT_m= S(m + \Frac 12 ) $ has crossing number $m+2$.
We have seen (example \ref{twh4}) that the only twist knots that are harmonic for $a=4$ are the
trefoil and the $\overline{5}_2$ knot.
The knot $\cT _m$  is not
harmonic for $a=3$ because $2^2 \not \equiv \pm 1 \Mod{2m+1}$ except when
$m=2$ (the figure-eight knot) or $m=1$ (trefoil).
From example \ref{twk}, we know that:
\bi
\item $\cT_{2k+1}$ can be parametrized by
$x=T_3(t), \, y=T_{3k+4}, \, z=C(t)$ where $\deg(C)=3k+5$.
\item $\cT_{2k}$ can be parametrized by
$x=T_3(t), \, y=T_{3k+2}, \, z=C(t)$ where $\deg(C)=3k+4.$
\ei
Using results of \ref{twk4}, we deduce other parametrizations
\bi
\item $\cT_{2k+1}$ can be represented by
$x=T_4(t), \  y= T_{2k+3} (t),\ z= C(t)$ where $C(t)$ is an odd polynomial
of degree $2k+5.$
\item $\cT_{2k}$ can be represented by
$ x= T_4(t), y= T_{2k+5} (t), z= C(t)$ where $C(t)$ is an odd polynomial of degree
$2k+7$.
\ei
\Pf
The proof is very similar to the preceding one,
except that there is an islet in the continued fractions for $2k$ even.
\EPf
\pn
Note that Chebyshev diagrams  we obtain ($a=3$ or $a=4$) for $\cT_{2k+1}$ have the same number
of crossing points: $3k+3$.
\begin{example}[The figure-eight knot]
$\cT_2$ is the figure-eight knot.
\begin{figure}[th]
\begin{center}
\begin{tabular}{ccc}
{\scalebox{.7}{\includegraphics{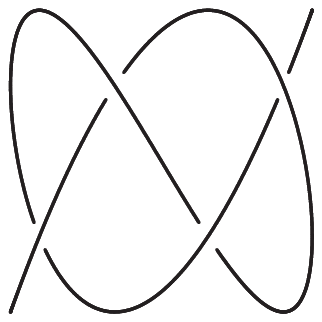}}}&
{\scalebox{.7}{\includegraphics{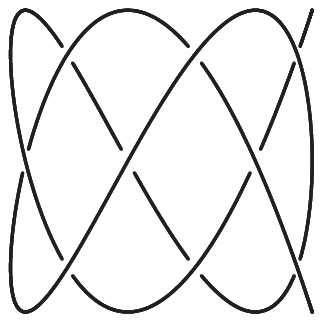}}}&
{\scalebox{.7}{\includegraphics{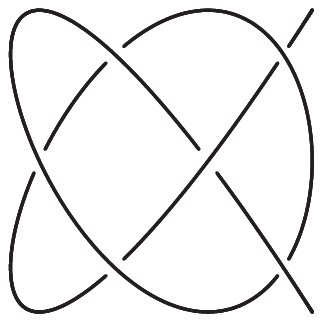}}}\\
$\overline{\H}(3,5,7)$&$S(\Frac 52)$&$S(-\Frac 58)$
\end{tabular}
\end{center}
\caption{The figure-eight knot}
\label{K41}
\end{figure}
Note that we obtain the figure-eight knot as the harmonic knot
$\H(3,5,7)$ or as a Chebyshev knot
$$x=T_4(t), \, y=T_7(t), z=
t \left( 10\,{t}^{2}-9 \right)  \left( 4\,{t}^{2}-3 \right)  \left( 3
\,{t}^{2}-2 \right)  \left( 2\,{t}^{2}-1 \right).
$$
But, considering $S(-5/8)$, we obtain a better parametrization
$$x=T_4(t), \, y = T_5(t), \, z= t \left( 11\,{t}^{2}-10 \right)  \left( 5\,{t}^{2}-4 \right)  \left( 5
\,{t}^{2}-1 \right).$$
\end{example}
\begin{example}[The 3-twist knot]
$\cT_3$ is the 3-twist knot $\overline{5}_2$. It is the harmonic knot $\H(4,5,7).$
It can also be parametrized by
$$x=T_3(t), \ y=T_7(t),\ z= t \left( 4\,t+3 \right)  \left( 3\,t+1 \right)  \left( 6\,t-5 \right)
 \left( 12\,{t}^{2}-11 \right)  \left( 2\,{t}^{2}-1 \right)$$
%or
%$x=T_4(t), \, y=T_5(t),\, z=t \left( 11\,{t}^{2}-10 \right)  \left( 2\,{t}^{2}-1 \right)  \left( 3
%\,{t}^{2}-1 \right)$
\begin{figure}[th]
\begin{center}
\begin{tabular}{cc}
{\scalebox{.7}{\includegraphics{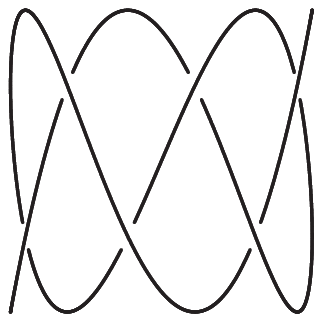}}}&
{\scalebox{.7}{\includegraphics{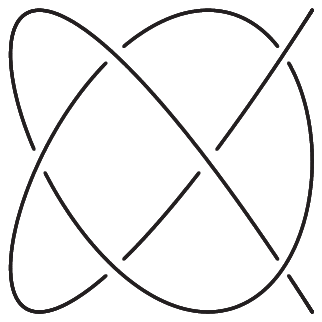}}}\\
$S(\Frac 72)$&$S(\Frac 74)$
\end{tabular}
\end{center}
\caption{Diagrams of the 3-twist knot $\overline{5}_2$}
\label{K52}
\end{figure}
\end{example}
\subsubsection*{Parametrizations of the generalized stevedore knots}
The stevedore knot $\cS_m= S (2m+2+ \Frac 1{2m} )$  can be represented by
$ x=T_3(t), \ y= T_{6m+2}(t), \ z = C(t) $ where $C(t)$ is a polynomial
of degree $6m+4$.

The stevedore knot $\cS_m= S (2m+2+ \Frac 1{2m} )$  can be represented by
$ x=T_4(t), \ y= T_{6m+3}(t), \ z = C(t) $ where $C(t)$ is an odd polynomial
of degree $ c= 10m +1.$
\pn
\Pf
The case $a=3$ is deduced from \ref{sk3} and Corollary \ref{degc3}.
The case $a=4$ is a consequence of Theorem \ref{gauss4}.
In this case $b=6m+3$, and the crossing number is $N= 4m+2.$
For $m=2k-1$ the number of islets in
$ \Frac { (4k-1)^2 }{4k} = (ASB)^{2k-1}(BSA)^k (\infty) $ is $\sigma = 2k-1=m.$
For $m= 2k,$ we also find $\sigma= m.$
Consequently we get
$ 3( 6m+3)+c -2 = 4 (4m+2)+ 12 m,$ that is,  $c= 10 m +1.$
The rest of the proof is analogous to the preceding ones.
\EPf
\begin{remark}
There is an algorithm to determine minimal Chebyshev diagrams for $a=3$ (Remark \ref{minb} and
Prop. \ref{ell}). When $a=4$, we can determine Chebyshev diagrams using Theorem \ref{cf1212} but
we do not know yet if they are minimal (consider for example $S(\Frac 9{14})$ and $S(-\Frac 58)$).
\end{remark}
\section{Proofs of theorems \ref{h3} and \ref{h4}}\label{proofs}
\subsection*{Proof of Theorem \ref{h3}}
\label{h3proof}
We study here the diagram of $\H(3,b,c)$ where $b=3n+1$ and
$c=2b-3\lambda$. The crossing points of the plane projection of
$\H(3,b,c)$ are obtained for pairs of values $(t,s)$
where
$t= \cos \bigl( \Frac m{3b} \pi \bigr), \ s = \cos \bigl( \Frac {m'}{3b} \pi \bigr).$
For $k = 0, \ldots, n-1$, let us consider
\bi
\item $A_{k}$  obtained for $m=3 k +1, \  m'= 2b-m.$
\item $B_{k}$  obtained for $ m =   3 k +2, \ m'= 2b + m$.
\item $C_{k}$  obtained for $ m = 2b - 3k - 3, \  m'= 4b-m$.
\ei
Then we have
\bi
\item $x(A_{k}) = \cos \bigl( \Frac {3 k +1}b \pi \bigr)$, \
$y(A_k) = \frac 12 (-1)^k$.
\item $x(B_{k} )= \cos \bigl( \Frac {3 k +2}b \pi \bigr)$, \
$y(B_k) = \frac 12 (-1)^{k+1}$.
\item
$x(C_{k}) = \cos \bigl( \Frac {3 k +3}b \pi \bigr)$, \
$y(C_k) = \frac 12 (-1)^{k}$.
\ei
\psfrag{a0}{\small $A_0$}\psfrag{b0}{\small $B_0$}\psfrag{c0}{\small $C_0$}%
\psfrag{a1}{\small $A_1$}%
\psfrag{an}{\small $A_{n-1}$}\psfrag{bn}{\small $B_{n-1}$}\psfrag{cn}{\small $C_{n-1}$}%
\begin{figure}[th]
\begin{center}
{\scalebox{.7}{\includegraphics{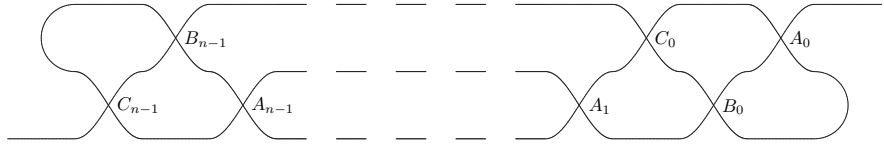}}}
\end{center}
\caption{$\H(3,3n+1,c)$, $n$ even}
\label{dh3}
\end{figure}
Hence our $3n$ points satisfy
$$x(A_{k-1}) > x(B_{k-1})
>x(C_{k-1}) >  x( A_{k}) > x(B_{k}) > x ( C_{k} ),\ k=1, \ldots, n-1.$$
Using the identity $T'_a( \cos \tau ) =  a  \Frac{ \sin a \tau}{\sin \tau },$
we get
$
x'(t)y'(t)  \sim \sin \bigl( \Frac mb \pi \bigr)  \sin \bigl( \Frac m3 \pi \bigr).
$
We obtain
\bi
\item[]
for $A_k$:
$
\begin{array}[t]{rcl}
{x'(t)y'(t)}&\sim& {\sin ( \Frac {3 k +1}b \pi ) \sin ( \Frac {3 k
+1}3 \pi )}  \sim (-1)^{k}.
\end{array}$
\item[] for $B_k$:
$
\begin{array}[t]{rcl}
x'(t)y'(t)&\sim& {\sin \bigl( \Frac {3 k +2}b \pi \bigr) \sin
\bigl( \Frac{3k+2}3 \pi \bigr)}
 \sim (-1)^{k}.
\end{array}$
\item[] for $C_k$:
$
\begin{array}[t]{rcl}
{x'(t)y'(t)}&\sim&
{\sin ( \Frac {2b -3 k -3}b \pi ) \sin (\Frac {2b-3
k -3}3 \pi )}\\
&\sim&
- {\sin (\Frac{3k+3}b \pi ) \sin ( -\Frac {3k+1}3  \pi ) } \sim
(-1)^{k}.
\end{array}$
\ei
The following identity will be useful in computing the sign of $z(t)-z(s).$
$$
T_c(t) - T_c(s) = 2 \sin \Bigl( \Frac{c}{6b}(m'-m) \pi  \Bigr)
\sin \Bigl( \Frac{c}{6b}(m+m') \pi   \Bigr).
$$
We have, with $c=2b -3\lambda$, $\theta= \Frac{\lambda}b \pi,$ (and $b=3n+1$ ),
\bi
\item[]
for $A_k$:
$ z(t)-z(s) =
-2 \sin c \Frac{\pi}3  \sin \Bigl( c \Frac{m-b}{3b} \pi \Bigr).$
But
\[
\sin  c \Frac {\pi}3 = \sin \Bigl( \Frac {6n+2-3 \lambda}3 \pi \Bigr)
= (-1)^{\lambda}  \sin \Frac {2 \pi }3,\label{c}
\]
and
$$ \sin \Bigl( c \Frac{b-m}{3b} \pi \Bigr)=
\sin \Bigl( (2- \Frac {3 \lambda}b )  \,  \Frac {b-m}3 \pi \Bigr) =
\sin \bigl( \Frac {\lambda}b (m-b) \pi \bigr) =
(-1)^{\lambda} \sin ( 3k+1) \theta.$$
We deduce that $z(t) - z(s) \sim \sin (3k+1) \theta$.
Finally, we obtain
$$ \sign{D(A_k)}= (-1)^k\sign{\sin (3k+1) \theta}.$$
\item
for $B_k$:
$
z(t)-z(s) =
2 \sin  c \Frac{\pi}{3}  \sin \bigl( \Frac cb \,. \,  \Frac{b+m}{3} \pi\bigr).
$
We have
\[
\sin \bigl( \Frac cb \cdot \Frac {b+m}3 \pi \bigr) &=&
\sin \bigl( ( 2-\Frac {3\lambda}b ) \Frac {b+m}3 \pi \bigr) \nonumber
\\
&=&
-\sin \bigl( \Frac {\lambda}b (b+m) \pi \bigr) =
(-1)^{\lambda+1} \sin  (3k+2) \theta. \nonumber
\]
Then, using Equation \ref{c}, we get $z(t) - z(s) \sim - \sin( 3k+2) \theta$,
and finally
$$
\sign{D(B_k)} = (-1)^{k+1} \sign{\sin(3k+2) \theta}.
$$
\item
for $C_k$:
$
\begin{array}[t]{rcl}
z(t)-z(s) &\sim&
\sin  \Frac {2 c }3 \pi  \sin \bigl( \Frac cb ( k+1)  \pi \bigr)\\
&\sim&  \sin   \Frac {4 \pi } 3  \sin \bigl(
(2- \Frac{3\lambda}b )( k+1) \pi \bigr) \sim     \sin (3k+3) \theta.
\end{array}$\\
We obtain
$$
\sign{ D(C_k)} = (-1)^k \sign{ \sin (3k+3) \theta }.
$$
\ei
These results give the Conway normal form.
If $n$ is odd, the Conway's signs of our points are
$$
\begin{array}{rcccl}
s(A_k) &\sim& (-1)^k D( A_k ) &\sim& { \sin (3k+1) \theta },\\
s(B_k)&\sim& (-1)^{k+1} D(B_k) &\sim& { \sin (3k+2) \theta }, \\
s(C_k) &\sim& (-1)^k D(C_k) &\sim& { \sin (3k+3) \theta }.
\end{array}
$$
In this case our result follows,  since the fractions
$ [e_1, e_2, \ldots, a_{3n} ] $ and  $(-1)^{3n+1} [a_{3n}, \ldots, a_1] $ define
the same knot.
If $n$ is even, the Conway's signs are
the opposite signs, and we also get the Schubert fraction of our knot.
\pn
Since $ 0<\theta < \Frac {\pi}2 $, we see that
there are not two consecutive sign changes in our sequence.
We also see that the  first two terms are of the same sign, and so are the last two terms.
The Conway normal form is biregular and the total number of sign changes in this sequence is
$\lambda -1 $: the crossing number of our knot is then
$b- \lambda.$ Finally, we get $ \beta ^2 \equiv \pm 1 $ by Proposition \ref{palin}.
\EPf
\subsection*{Proof of Theorem \ref{h4}}
\label{h4proof}
The crossing points of the plane projection of
$\H= \H (4,b,c)$ are obtained for parameter pairs $ ( t, s) $ where
$ t=\cos \bigl( \Frac m{4b} \pi \bigr), \ s= \cos \bigl( \Frac {m'}{4b} \pi \bigr).$
We shall denote
$ \lambda = \Frac {3b-c}4,$
( or $ c= 3b - 4 \lambda$) and $ \theta = \Frac {\lambda}b \pi.$
We will consider the two following  cases.
\pn
{\bf The case $\mathbf{b = 4n+1}$.}
\pn
For $k=0, \ldots, n-1$,  let us consider the following crossing points
\bi
\item
$A_k$ corresponding to $ m=4k+1, \  m'= 2b-m,$
\item
$B_k$ corresponding to $ m= 4k+2, \  m'= 4b-m,$
\item
$C_k$ corresponding to $ m= 4k+3, \ m' = 2b+m,$
\item
$ D_k$ corresponding to $ m = 2b - 4(k+1), \  m' = 4b-m.$
\ei
Then we have
\bi
\item
$ x( A_k)= \cos \bigl( \Frac {4k+1}b \pi \bigr), \  \  y(A_k) = (-1)^k \cos \Frac \pi 4 \neq 0, $
\item
$ x(B_k) = \cos \bigl( \Frac {4k+2}b \pi \bigr), \   \  y(B_k)=0,$
\item
$ x( C_k) = \cos \bigl( \Frac {4k+3}b \pi \bigr), \  \
 y(C_k) = (-1)^k \cos \Frac {3 \pi }4 \neq 0, $
\item
$x(D_k) = \cos \bigl( \Frac {4k+4}b \pi \bigr), \  \  y(D_k)=0.$
\ei
\psfrag{a0}{\small $A'_0$}\psfrag{b0}{\small $B_0$}\psfrag{c0}{\small $C_0$}\psfrag{d0}{\small $D_0$}%
\psfrag{an}{\small $A_{n-1}$}\psfrag{bn}{\small $B_{n-1}$}\psfrag{cn}{\small $C'_{n-1}$}\psfrag{dn}{\small $D_{n-1}$}%
\psfrag{aa0}{\small $A_0$}\psfrag{cc0}{\small $C'_0$}%
\psfrag{aan}{\small $A'_{n-1}$}\psfrag{ccn}{\small $C_{n-1}$}%
\begin{figure}[th]
\begin{center}
{\scalebox{.7}{\includegraphics{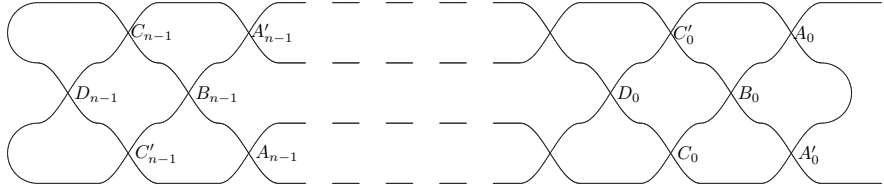}}}
\end{center}
\caption{$\H(4,4n+1,c)$, $n$ even}
\label{dh4}
\end{figure}
Hence our $4n$ points satisfy
\$x(A_{k-1}) > x(B_{k-1}) > x( C_{k-1}) >  x( D_{k-1})>&&\\
\quad\quad x(A_k) > x( B_k) > x(C_k ) > x( D_k),&& k=1, \ldots, n-1.
\$
We remark that these points together with the symmetric points $A'_k$ (resp. $C'_k$) of $A_k$
(resp. $C_k$) with respect to the $y-$axis form the totality of the crossing points.

The Conway sign of  a crossing point  $M$ is
$s(M)= \sign{D(M)}$ if $y(M)=0,$ and
$s(M)=- \sign{D(M)}$ if $y(M) \neq 0.$

By symmetry, we have $s(A'_k) = s(A_k)$ and  $s(C'_k) = s(C_k)$ because symmetric
points correspond to opposite parameters.
The Conway  form of $\H$ is then (see paragraph \ref{cf}) :
$$
C\Bigl( s (D_{n-1}), 2s(C_{n-1}), s(B_{n-1}),2s( A_{n-1}), \ldots, s (B_0), 2s (A_0) \Bigr).
$$
Using the identity $ T' _a ( \cos \tau ) = a \Frac { \sin a \tau }{ \sin \tau},$
we get
$ x'(t) y'(t)\sim \sin \bigl( \Frac mb \pi \bigr)  \sin \bigl( \Frac {m}4 \pi  \bigr).$
Consequently,
\bi
\item
For $A_k$ we have
$ x'(t) y'(t) \sim
\sin \bigl( \Frac {4k+1}b \pi \bigr) \sin \bigl( \Frac {4k+1}4 \pi \bigl)
\sim(-1)^k.$
\item
Similarly, for $B_k$ and $C_k$ we get
$ x'(t) y'(t) \sim (-1)^k. $
\item
For $D_k$ we get
$ x'(t) y'(t)\sim
\sin \bigl( \Frac {2b-4k-4}b \pi \bigr)
\sin \bigl( \Frac {2b-4k-4}4 \pi \bigr)\sim (-1)^{k+1}.$
\ei
On the other hand, at the crossing points we have
$$ z(t) - z(s) = 2 \sin \Bigl( \Frac c{8b} (m'-m) \pi \Bigr) \,
\sin  \Bigl( \Frac c{8b} (m+m') \pi  \Bigr).
$$
We obtain the signs of our crossing points,
with $c = 3b-4 \lambda, \  \theta= \Frac {\lambda}b,t.$
\bi
\item
For $A_k$ we get:
$ z(t )-z(s)= 2 \sin \Frac cb (n-k) \pi\, \sin  c \Frac {\pi}4.$\\
We have
$ \sin c \Frac {\pi}4 = \sin  \Frac { 12 n +3 - 4 \lambda}4  \pi =
(-1)^{n+ \lambda} \sin  \Frac {3 \pi }4 \sim (-1)^{n+ \lambda}$\\
and also
$
\begin{array}[t]{rcl}
\sin \bigl( \Frac cb (n-k) \pi  \bigr) &=&
\sin \Bigl(  \bigl( 3 - \Frac {4 \lambda}b \bigr) \bigl( n-k) \pi \Bigr)\\
&=& (-1)^{n+k} \sin \Bigl(  \Frac {4k-4n}b \lambda \pi  \Bigr) = (-1) ^{ n+k+ \lambda} \sin ( 4k+1) \theta
\end{array}$.\\
Consequently, the sign of $A_k$ is
$$ s(A_k) = - \sign{\sin ( 4 k+1) \theta }.$$
\item
For $B_k$, we have:
$z(t) - z(s) = 2 \sin \bigl(  \Frac cb ( 2n-k) \pi \bigr)
\sin c \Frac {\pi}2 = - 2 \sin \bigl(  \Frac cb ( 2n-k) \pi \bigr).$\\
But
$
\begin{array}[t]{rcl}
\sin \bigl( \Frac cb ( 2n-k) \pi \bigr) &=&
\sin \Bigl( \bigl( 3 - \Frac{4 \lambda}b \bigr) \bigl( 2n-k \bigr) \pi \Bigr) \\
&=& (-1)^{k} \sin \Bigl( \Frac {\lambda}b ( 4k-8n) \pi \Bigr) = (-1)^k \sin (4k+2) \theta.
\end{array}$.\\
Therefore the sign of $B_k$ is
$$
s (B_k) = - \sign{\sin (4k+2) \theta}.$$
\item For $C_k$:
$ z(t) - z(s) = 2 \sin \bigl( \Frac c4 \pi \bigr)
\sin \bigl(  \Frac cb (n+k+1) \pi  \bigr). $\\
We know that $ \sin \Frac {c \pi} 4 \sim (-1)^{n+ \lambda}$. Let us compute the second factor
$$
\begin{array}{rcl}
\sin \Bigl(  \bigl( 3 - \Frac {4 \lambda}b \bigr) \bigl(  n+k+1 \bigr) \pi\Bigr) &=&
(-1)^{n+k} \sin \Bigl( \Frac {\lambda}b \bigl( 4n+4k+4 \bigr) \pi \Bigr) \\
&=&(-1)^{n+k} \sin \Bigl( \Frac {\lambda}b ( b+4k+3) \pi \Bigr) \\
&=& ( -1)^{n+k+ \lambda} \sin (4k+3)  \theta.
\end{array}$$
Hence
$$ s( C_k) = - \sign{\sin( 4k+3) \theta}.$$
\item
For $D_k$:
$
\begin{array}[t]{rcl}
z(t) -z(s)&=& 2 \sin \bigl( \Frac cb (k+1) \pi \bigr) \sin \bigl( c \Frac {\pi}2 \bigr)\\
&=& 2 \sin \Bigl(  \bigl( 3- \Frac {4 \lambda }b   \bigr) \bigl( k+1 \bigr) \pi \Bigr) =
(-1)^{k} \sin ( 4k+4) \theta .
\end{array}$.\\
We conclude
$$s (D_k)=
- \sign{\sin( 4k+1) \theta}.$$
\ei
This completes the computation of our Conway normal form of $H$ in this first case.
\pn
{\bf The case $\mathbf{b = 4n+3}$.}\\
Here, the diagram is different.
Let us consider the following $4n+2$ crossing points.

For $k= 0, \ldots, n$
\bi
\item
$ A_k$ corresponding to $ m=4k+1, \  m' = 2b+m,$
\item
$B_k$ corresponding to $ m=4k+2, \  m' = 4b-m.$
\ei
For $k=0, \ldots, n-1$
\bi
\item
$ C_k$ corresponding to $ m=4k+3, \ m' =2b-m,$
\item
$ D_k$ corresponding to $ m = 2b + 4 (k+1),  \  m'=4b-m.$
\ei
These points are chosen so that
$$
x(A_0) > x(B_0)> x(C_0) > x(D_0) > \cdots > x(D_{n-1}) > x( A_n) > x(B_n),
$$
and we have $\sign{x'(t) y'(t)}= (-1)^k.$
\bi
\item For $A_k$ we get
$$
z(t)-z(s)= 2 \sin \bigl( c \Frac {\pi}4 \bigr)
\sin \bigl(  \Frac cb ( n+k+1) \pi \bigr).
$$
We  easily get $ \sign{\sin c \Frac {\pi}4 } = (-1)^{n+ \lambda}.$
We also get
$$
\begin{array}{rcl}
\sin \bigl( \Frac cb ( n+k+1) \pi \bigr)&=&
\sin \Bigl( \bigl( 3 - \Frac {4 \lambda}b \bigr) \bigl( n+k+1 \bigr) \pi \Bigr)\\
&=&
(-1)^{n+k} \sin \bigl( \Frac {\lambda}b ( b+4k+1) \pi \bigr) =
(-1)^{n+k+\lambda} \sin ( 4k+1) \theta.
\end{array}
$$
Hence the sign of $A_k$ is
$$
s(A_k)= - \sign{\sin  (4k+1) \theta}.
$$
\item For $B_k$ we get
$$
z(t)-z(s)= 2 \sin \bigl(  \Frac cb ( 2n+1-k) \pi \bigr)
\sin c \Frac {\pi}2.
$$
We have $ \sin \bigl( c \Frac {\pi}2 \bigr) =1 >0, $ and
$$
\begin{array}{rcl}
\sin \bigl( \Frac cb (2n+1-k) \pi \bigr)&=&
\sin \Bigl( \bigl( 3- \Frac { 4 \lambda}b \bigr) \bigl( 2n+1 -k \bigr) \pi \Bigr)\\
&=&(-1)^{k+1} \sin \bigl( \Frac {\lambda}b ( 4k-8n-4) \pi \bigr) =
(-1)^{k+1} \sin (4k+2) \theta.
\end{array}
$$
Then, the sign of $B_k$ is
$$
s(B_k) = - \sign{\sin(4k+2) \theta}.
$$
\item For $  C_k$ we have
$$
z(t) - z(s) = 2 \sin \bigl( \Frac cb (n-k) \pi  \bigr) \sin c \Frac {\pi}4.
$$
We get
$$
\begin{array}{rcl}
\sin \bigl( \Frac cb (n-k) \pi \bigr) &=&
\sin \Bigl( \bigl( 3- \Frac {4 \lambda}b \bigr) \bigl( n-k \bigr) \pi \Bigr) \\
&=& (-1)^{n+k} \sin \bigl( \Frac {4k-4n}b \lambda \pi \bigr)=
(-1)^{n+k+ \lambda} \sin (4k+3) \theta.
\end{array}
$$
The sign of $C_k$ is then
$$
s( C_k) =
- \sign{\sin ( 4k+3) \theta}.
$$
\item For $D_k$ we get
$$
z(t)-z(s)= 2 \sin \bigl( -\Frac{c}b (k+1) \pi \bigr) \sin c \Frac {\pi}2.
$$
We have $ \sin c \Frac {\pi}2 >0.$
We also have
$$
\begin{array}{rcl}
\sin \bigl( -\Frac cb (k+1) \pi ) &=&
\sin \Bigl( \bigl( \Frac {4 \lambda}b -3 \bigr) \bigl( k+1 \bigr) \pi  \Bigr)
(-1)^{k+1} \sin(4k+4) \theta.
\end{array}
$$
Consequently, the sign of $D_k$ is
$$
s(D_k)= - \sign{\sin(4k+4) \theta }.
$$
\ei
This concludes the computation of the  Conway normal form of $\H(4,b,c).$
\pn
If $ b<c< 3b,$ we get
$ \lambda < \Frac b2,$ and then  $ \theta < \Frac {\pi} 2.$
Consequently, our sequence is biregular.
Furthermore, the total number of sign changes is $ \lambda-1.$
We conclude that the crossing number is
$ N= \Frac { 3 (b-1) }2 - (\lambda -1)= \Frac {3b+c-2}4.$
The fact that $ \beta ^2 \equiv \pm2 \Mod \alpha $ is a consequence of
Proposition \ref{palin4}.
\EPf
%\section{Conclusion}

\vfill
\pn
\hrule width 5cm height 2pt
\pn
Pierre-Vincent Koseleff, \\
\'Equipe-project INRIA Salsa \& Universit{\'e} Pierre et Marie Curie (UPMC-Paris 6)\\
e-mail: {\tt koseleff@math.jussieu.fr}
\pn
Daniel Pecker, \\
Universit{\'e} Pierre et Marie Curie (UPMC-Paris 6)\\
e-mail: {\tt pecker@math.jussieu.fr}

\begin{thebibliography}{ZZ99}
%\bibitem [Bi]{Bi}
%J. S. Birman, {\it Braids, Links and Mapping Class Groups,} Ann.
%of Math.Studies 82, Princeton University Press, 1974.
%
\bibitem [BZ]{BZ}
G. Burde, H. Zieschang,
{\it Knots}, Walter de Gruyter, 2003
%
\bibitem [BDHZ] {BDHZ}
A. Boocher, J. Daigle, J. Hoste, W. Zheng,
{\it Sampling Lissajous and Fourier knots},
arXiv:0707.4210, (2007).
%
\bibitem[BHJS]{BHJS}
M. G. V. Bogle, J. E. Hearst, V. F. R. Jones, L. Stoilov, {\it
Lissajous knots}, Journal of Knot Theory and its Ramifications,
3(2): 121-140, (1994).
%
\bibitem [Com]{Com}
E. H. Comstock,
{\it The Real Singularities of Harmonic Curves of three Frequencies},
Trans. of the Wisconsin Academy of Sciences, Vol XI : 452-464, (1897).
%
\bibitem [Con]{Con}
J. H. Conway,
{\it An enumeration of knots and links, and some of their algebraic properties},
Computational Problems in Abstract Algebra (Proc. Conf., Oxford, 1967), 329--358 Pergamon,
Oxford (1970)
%
\bibitem [Cr]{Cr}
P. R. Cromwell, {\it Knots and links},
Cambridge University Press, Cambridge, 2004. xviii+328 pp.
%
%\bibitem [ES]{ES}
%C. Ernst, D. W. Sumners, {\it The Growth of the Number of Prime Knots}, Math. Proc.
%Cambridge Philos. Soc. {\bf 102},
%303-315, (1987).
%
\bibitem[Fi]{Fi}
G. Fischer, {\it Plane Algebraic Curves}, A.M.S. Student Mathematical
Library Vol 15, 2001.
%
\bibitem[FF]{FF}
G. Freudenburg, J. Freudenburg,
{\it Curves defined by Chebyshev polynomials}, 19 p., (2009),
\verb+arXiv:0902.3440+
%
\bibitem[HK]{HK}
R. Hartley, A. Kawauchi, {\it Polynomials of amphicheiral knots},
Math. Ann. {\bf 243 (1)}: 63-70 (1979)
%
%\bibitem[HW]{HW}
%G.H. Hardy, E.M. Wright, {\it An Introduction to the Theory of Numbers,}
%4th edition, Oxford University Press, 1960.
%
\bibitem[HZ]{HZ}
J. Hoste, L. Zirbel, {\it Lissajous knots and knots with Lissajous
  projections,}
\verb+arXiv:math.GT/0605632v1+, (2006).
To appear in Kobe Journal of mathematics, vol 24, n$^{\rm o}2$
%
\bibitem[JP]{JP}
V. F. R. Jones, J. Przytycki, {\it Lissajous  knots and billiard knots,}
Banach Center Publications, 42:145-163, (1998).
%
\bibitem[Kaw]{Kaw}
A. Kawauchi, editor, {\it A Survey of Knot Theory}, Birh{\"a}user, 1996.
%
\bibitem[KP1]{KP1}
P. -V. Koseleff, D. Pecker, {\it On polynomial Torus Knots},
Journal of Knot Theory and its Ramifications, Vol. {\bf 17 (12)}
(2008), 1525-1537
%
\bibitem[KP2]{KP2}
P. -V. Koseleff, D. Pecker, {\it A construction of polynomial torus
 knots}, to appear in Journal of AAECC\\
\verb+http://arxiv.org/abs/0712.2408+
%
\bibitem[KP3]{KP3}
P. -V. Koseleff, D. Pecker, {\it Chebyshev knots},
\verb+arXiv:0812.1089+, (2008).
%
\bibitem[KPR]{KPR}
P. -V. Koseleff, D. Pecker, F. Rouillier, {\it The first rational Chebyshev knots}, Conference MEGA 2009, Barcelona.
%
\bibitem[La1]{La1}
C. Lamm, {\it There are infinitely many Lissajous knots,}
Manuscripta Math., 93: 29-37, (1997).
%
\bibitem[La2]{La2}
C. Lamm, {\it Fourier Knots,} preprint
%
%\bibitem[Li]{Li}
%Livingston, C., {\it Knot Theory}, Washington, DC: Math. Assoc. Amer., 1993.
%
\bibitem[Mi]{Mi}
R. Mishra, {\it Polynomial representations of strongly-invertible knots and
strongly-negative amphicheiral knots, }  Osaka J. Math., 43 625-639 (2006).

\bibitem[Mu]{Mu}
K. Murasugi, {\it Knot Theory and its Applications}, Boston,
Birkh{\"a}user, 341p., 1996.
%
\bibitem[P1]{P1}
D. Pecker, {\it Simple constructions of algebraic curves with nodes},
Compositio Math.  87  (1993),  no. 1, 1--4.
%
\bibitem [P2]{P2}
D. Pecker, {\it Sur le genre arithm{\'e}tique des courbes rationnelles. (French)
[On the arithmetic genus of rational curves]},
Ann. Inst. Fourier (Grenoble)  46  (1996),  no. 2, 293--306.
%
\bibitem[RS]{RS}
A. Ranjan and R. Shukla, {\it On polynomial representation of
torus knots,} Journal of knot theory and its ramifications, Vol. 5
(2) (1996) 279-294.
%
\bibitem[Ro]{Ro}
D. Rolfsen, {\it Knots and Links},
Math. Lecture Series 7, Publish or Perish, 1976.
%
\bibitem[Sh]{Sh}
A.R. Shastri, {\it Polynomial representation of knots},
 T{\^o}hoku Math. J. {\bf 44 } (1992), 11-17.
%
\bibitem[St]{St}
A. Stoimenow, {\it Generating functions, Fibonacci numbers and rational knots},
J. Algebra {\bf 310(2)} (2007), 491--525.
%
\bibitem[Tu]{Tu}
J.C. Turner, {\it On a class of knots with Fibonacci invariant numbers},
Fibonacci Quart. {\bf 24} (1986), n$^{\rm o}1$,  61-66.
%
\bibitem[Va]{Va}
V. A. Vassiliev, {\it Cohomology of knot spaces}, Theory of
singularities and its Applications, Advances Soviet Maths Vol 1
(1990)
\end{thebibliography}
\end{document}